\documentclass[12pt]{article}

\usepackage[pdfa]{hyperref}

\usepackage[utf8]{inputenc}
\usepackage{eurosym}
\usepackage{latexsym,amssymb,amsmath}
\usepackage{tikz}
\usepackage{amsmath,amsfonts,latexsym, amssymb}
\usepackage{mathrsfs}
\usepackage{graphics}
\usepackage{float}
\usepackage{stmaryrd}
\usepackage{geometry}
\usepackage{lscape}
\usepackage{comment}
\usepackage{listings}
\newtheorem{thm}{Theorem}[section]
\newtheorem{prop}[thm]{Proposition}
\newtheorem{lemma}[thm]{Lemma}
\newtheorem{cor}[thm]{Corollary}
\newtheorem{conj}[thm]{Conjecture}

\newtheorem{claim}[thm]{Claim}

\newcommand{\F}{\mathcal{F}}

\newcommand{\R}{\mathbb{R}}
\newcommand{\im}{\text{im}}
\newcommand{\ext}{\text{ex}}
\newcommand{\Hom}{\text{Hom}}
\newcommand{\OPT}{\text{OPT}}
\newcommand{\fl}{\text{fl}}
\newcommand{\cl}{\text{cn}}

\newenvironment{proof}[1][Proof]{\textbf{#1.} }{\ \rule{0.5em}{0.5em}}

\newcommand{\vc}[1]{\ensuremath{\vcenter{\hbox{#1}}}}
\newcommand{\fourone}{\vc{\includegraphics[page=1, scale = 0.75]{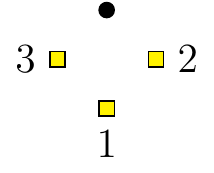}}}
\newcommand{\fourtwo}{\vc{\includegraphics[page=2, scale = 0.75]{figs4label.pdf}}}
\newcommand{\fourthree}{\vc{\includegraphics[page=3, scale = 0.75]{figs4label.pdf}}}
\newcommand{\fourfour}{\vc{\includegraphics[page=4, scale = 0.75]{figs4label.pdf}}}
\newcommand{\fourfive}{\vc{\includegraphics[page=5, scale = 0.75]{figs4label.pdf}}}
\newcommand{\foursix}{\vc{\includegraphics[page=6, scale = 0.75]{figs4label.pdf}}}
\newcommand{\fourseven}{\vc{\includegraphics[page=7, scale = 0.75]{figs4label.pdf}}}
\newcommand{\foureight}{\vc{\includegraphics[page=8, scale = 0.75]{figs4label.pdf}}}
\newcommand{\fourforkthree}{\vc{\includegraphics[page=15, scale = 0.75]{figs4label.pdf}}}
\newcommand{\diamondforkthree}{\vc{\includegraphics[page=16, scale = 0.75]{figs4label.pdf}}}

\newcommand{\Fa}{\vc{\includegraphics[page=1, scale = 0.4]{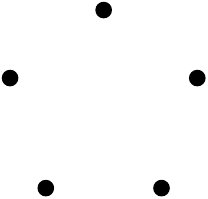}}}
\newcommand{\Fb}{\vc{\includegraphics[page=2, scale = 0.4]{figs5unlabeled.pdf}}}
\newcommand{\Fc}{\vc{\includegraphics[page=3, scale = 0.4]{figs5unlabeled.pdf}}}
\newcommand{\Fd}{\vc{\includegraphics[page=4, scale = 0.4]{figs5unlabeled.pdf}}}
\newcommand{\Fe}{\vc{\includegraphics[page=5, scale = 0.4]{figs5unlabeled.pdf}}}
\newcommand{\Ff}{\vc{\includegraphics[page=6, scale = 0.4]{figs5unlabeled.pdf}}}
\newcommand{\Fg}{\vc{\includegraphics[page=7, scale = 0.4]{figs5unlabeled.pdf}}}
\newcommand{\Fh}{\vc{\includegraphics[page=8, scale = 0.4]{figs5unlabeled.pdf}}}
\newcommand{\Fi}{\vc{\includegraphics[page=9, scale = 0.4]{figs5unlabeled.pdf}}}
\newcommand{\Fj}{\vc{\includegraphics[page=10, scale = 0.4]{figs5unlabeled.pdf}}}
\newcommand{\Fk}{\vc{\includegraphics[page=11, scale = 0.4]{figs5unlabeled.pdf}}}
\newcommand{\Fl}{\vc{\includegraphics[page=12, scale = 0.4]{figs5unlabeled.pdf}}}
\newcommand{\Fm}{\vc{\includegraphics[page=13, scale = 0.4]{figs5unlabeled.pdf}}}
\newcommand{\Fn}{\vc{\includegraphics[page=14, scale = 0.4]{figs5unlabeled.pdf}}}
\newcommand{\Fo}{\vc{\includegraphics[page=15, scale = 0.4]{figs5unlabeled.pdf}}}
\newcommand{\Fp}{\vc{\includegraphics[page=16, scale = 0.4]{figs5unlabeled.pdf}}}
\newcommand{\Fq}{\vc{\includegraphics[page=17, scale = 0.4]{figs5unlabeled.pdf}}}
\newcommand{\Fr}{\vc{\includegraphics[page=18, scale = 0.4]{figs5unlabeled.pdf}}}
\newcommand{\Fs}{\vc{\includegraphics[page=19, scale = 0.4]{figs5unlabeled.pdf}}}
\newcommand{\Ft}{\vc{\includegraphics[page=20, scale = 0.4]{figs5unlabeled.pdf}}}
\newcommand{\Fu}{\vc{\includegraphics[page=21, scale = 0.4]{figs5unlabeled.pdf}}}
\newcommand{\Fv}{\vc{\includegraphics[page=22, scale = 0.4]{figs5unlabeled.pdf}}}
\newcommand{\Fw}{\vc{\includegraphics[page=23, scale = 0.4]{figs5unlabeled.pdf}}}
\newcommand{\Fx}{\vc{\includegraphics[page=24, scale = 0.4]{figs5unlabeled.pdf}}}
\newcommand{\Fy}{\vc{\includegraphics[page=25, scale = 0.4]{figs5unlabeled.pdf}}}
\newcommand{\Fz}{\vc{\includegraphics[page=26, scale = 0.4]{figs5unlabeled.pdf}}}
\newcommand{\Faa}{\vc{\includegraphics[page=27, scale = 0.4]{figs5unlabeled.pdf}}}
\newcommand{\Fbb}{\vc{\includegraphics[page=28, scale = 0.4]{figs5unlabeled.pdf}}}
\newcommand{\Fcc}{\vc{\includegraphics[page=29, scale = 0.4]{figs5unlabeled.pdf}}}
\newcommand{\Fdd}{\vc{\includegraphics[page=30, scale = 0.4]{figs5unlabeled.pdf}}}
\newcommand{\Fee}{\vc{\includegraphics[page=31, scale = 0.4]{figs5unlabeled.pdf}}}
\newcommand{\Fff}{\vc{\includegraphics[page=32, scale = 0.4]{figs5unlabeled.pdf}}}
\newcommand{\Fgg}{\vc{\includegraphics[page=33, scale = 0.4]{figs5unlabeled.pdf}}}
\newcommand{\Fhh}{\vc{\includegraphics[page=34, scale = 0.4]{figs5unlabeled.pdf}}}

\newcommand{\Fasml}{\vc{\includegraphics[page=1, scale = 0.2]{figs5unlabeled.pdf}}}
\newcommand{\Fbsml}{\vc{\includegraphics[page=2, scale = 0.2]{figs5unlabeled.pdf}}}
\newcommand{\Fcsml}{\vc{\includegraphics[page=3, scale = 0.2]{figs5unlabeled.pdf}}}
\newcommand{\Fdsml}{\vc{\includegraphics[page=4, scale = 0.2]{figs5unlabeled.pdf}}}
\newcommand{\Fesml}{\vc{\includegraphics[page=5, scale = 0.2]{figs5unlabeled.pdf}}}
\newcommand{\Ffsml}{\vc{\includegraphics[page=6, scale = 0.2]{figs5unlabeled.pdf}}}
\newcommand{\Fgsml}{\vc{\includegraphics[page=7, scale = 0.2]{figs5unlabeled.pdf}}}
\newcommand{\Fhsml}{\vc{\includegraphics[page=8, scale = 0.2]{figs5unlabeled.pdf}}}
\newcommand{\Fisml}{\vc{\includegraphics[page=9, scale = 0.2]{figs5unlabeled.pdf}}}
\newcommand{\Fjsml}{\vc{\includegraphics[page=10, scale = 0.2]{figs5unlabeled.pdf}}}
\newcommand{\Fksml}{\vc{\includegraphics[page=11, scale = 0.2]{figs5unlabeled.pdf}}}
\newcommand{\Flsml}{\vc{\includegraphics[page=12, scale = 0.2]{figs5unlabeled.pdf}}}
\newcommand{\Fmsml}{\vc{\includegraphics[page=13, scale = 0.2]{figs5unlabeled.pdf}}}
\newcommand{\Fnsml}{\vc{\includegraphics[page=14, scale = 0.2]{figs5unlabeled.pdf}}}
\newcommand{\Fosml}{\vc{\includegraphics[page=15, scale = 0.2]{figs5unlabeled.pdf}}}
\newcommand{\Fpsml}{\vc{\includegraphics[page=16, scale = 0.2]{figs5unlabeled.pdf}}}
\newcommand{\Fqsml}{\vc{\includegraphics[page=17, scale = 0.2]{figs5unlabeled.pdf}}}
\newcommand{\Frsml}{\vc{\includegraphics[page=18, scale = 0.2]{figs5unlabeled.pdf}}}
\newcommand{\Fssml}{\vc{\includegraphics[page=19, scale = 0.2]{figs5unlabeled.pdf}}}
\newcommand{\Ftsml}{\vc{\includegraphics[page=20, scale = 0.2]{figs5unlabeled.pdf}}}
\newcommand{\Fusml}{\vc{\includegraphics[page=21, scale = 0.2]{figs5unlabeled.pdf}}}
\newcommand{\Fvsml}{\vc{\includegraphics[page=22, scale = 0.2]{figs5unlabeled.pdf}}}
\newcommand{\Fwsml}{\vc{\includegraphics[page=23, scale = 0.2]{figs5unlabeled.pdf}}}
\newcommand{\Fxsml}{\vc{\includegraphics[page=24, scale = 0.2]{figs5unlabeled.pdf}}}
\newcommand{\Fysml}{\vc{\includegraphics[page=25, scale = 0.2]{figs5unlabeled.pdf}}}
\newcommand{\Fzsml}{\vc{\includegraphics[page=26, scale = 0.2]{figs5unlabeled.pdf}}}
\newcommand{\Faasml}{\vc{\includegraphics[page=27, scale = 0.2]{figs5unlabeled.pdf}}}
\newcommand{\Fbbsml}{\vc{\includegraphics[page=28, scale = 0.2]{figs5unlabeled.pdf}}}
\newcommand{\Fccsml}{\vc{\includegraphics[page=29, scale = 0.2]{figs5unlabeled.pdf}}}
\newcommand{\Fddsml}{\vc{\includegraphics[page=30, scale = 0.2]{figs5unlabeled.pdf}}}
\newcommand{\Feesml}{\vc{\includegraphics[page=31, scale = 0.2]{figs5unlabeled.pdf}}}
\newcommand{\Fffsml}{\vc{\includegraphics[page=32, scale = 0.2]{figs5unlabeled.pdf}}}
\newcommand{\Fggsml}{\vc{\includegraphics[page=33, scale = 0.2]{figs5unlabeled.pdf}}}
\newcommand{\Fhhsml}{\vc{\includegraphics[page=34, scale = 0.2]{figs5unlabeled.pdf}}}

\newcommand{\threesigone}{\vc{\includegraphics[page=12, scale = 0.75]{figs4label.pdf}}}
\newcommand{\threesigtwo}{\vc{\includegraphics[page=13, scale = 0.75]{figs4label.pdf}}}
\newcommand{\threesigthree}{\vc{\includegraphics[page=14, scale = 0.75]{figs4label.pdf}}}

\newcommand{\figsplusone}{\vc{\includegraphics[page=1, scale = 1]{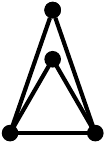}}}

\newcommand{\Dist}{\text{Dist}}

\title{Maximizing five-cycles in $K_r$-free graphs}
\author{Bernard Lidick\'y\thanks{Department of Mathematics, Iowa State University. Ames, IA, USA. E-mail: \texttt{lidicky@iastate.edu}. Supported in part by NSF grant DMS-1855653.}
\and 
Kyle Murphy\thanks{Department of Mathematics, Iowa State University. Ames, IA, USA. E-mail: \texttt{kylem2@iastate.edu}}
}
\date{\today}

\begin{document}

\maketitle

\begin{abstract}
The Pentagon Problem of Erd\H{o}s problem asks to find an $n$-vertex triangle-free graph that is maximizing the number of 5-cycles.
The problem was solved using flag algebras by Grzesik and independently by Hatami, Hladk\'{y}, Kr\'{a}l', Norin, and Razborov. 
Recently, Palmer suggested a more general problem of maximizing the number of 5-cycles in  $K_{k+1}$-free graphs. 
Using flag algebras, we show that every $K_{k+1}$-free graph of order $n$ contains at most
$$\frac{1}{10k^4}(k^4 - 5k^3 + 10k^2 - 10k + 4)n^5 + o(n^5)$$
copies of $C_5$ for any $k \geq 3$, with the Tur\'{a}n graph being the extremal graph for large enough  $n$.
\end{abstract}

\section{Introduction}
All graphs in this paper are simple.
Let $G$, $H$, and $F$ be graphs. We define $\nu(H,G)$ as the number of (possibly non-induced) subgraphs of $G$ isomorphic to $H$. If $G$ does not contain any subgraph isomorphic to $F$, then we say that $G$ is \emph{$F$-free}. Let $\ext(n,H,F)$ denote the maximum value of $\nu(H,G)$ among all $F$-free graphs $G$ on $n$ vertices. The function $\ext(n,H,F)$ is well-studied when $H$ is an edge. As such, it is convention when $H=K_2$ to let $\ext(n,F)$ denote $\ext(n,K_2,F)$. The value of $\ext(n,F)$ for any graph $F$ is called the \emph{Tur\'{a}n number} of $F$. 

One of the first results in extremal graph theory was Mantel's Theorem~\cite{mantel} which states that for all $n \geq 3$, $\ext(n,K_3) \leq \lfloor \frac{n^2}{4} \rfloor$. 
When $k \geq 3$, the value of $\ext(n,K_{k+1})$ was determined by Tur\'{a}n. 
\begin{thm}[Tur\'{a}n's Theorem~\cite{turan}]
For all $k \geq 3$, and all $n$,
\[ 
\ext(n,K_{k+1}) \leq \left(1 - \frac{1}{k}\right)\frac{n^2}{2}. 
\]
Moreover, the Tur\'{a}n graph $T_k(n)$, which is the complete balanced $(k-1)$-partite graph on $n$ vertices, is the unique $K_{k+1}$-free graph on $n$ vertices which contains the maximum possible number of edges. 
\end{thm}

The Erd\H{o}s-Stone-Simonovits Theorem~\cite{ErdosStone} determined the asymptotic value of $\ext(n,F)$ when $F$ is not a complete graph. Let $\chi(F)$ denote the chromatic number of $F$. Then for all $F$ for which $\chi(F) \geq 2$,

$$\ext(n,F) = \frac{\chi(F) - 2}{2\left(\chi(F) - 1\right)}n^2 + o(n^2).$$

The systematic study of the function $\ext(n,H,F)$ was initiated by  Alon and Shikhelman~\cite{Alon}, although there were some prior results. When $t < k+1$, Zykov~\cite{Zykov} showed that the Tur\'{a}n graph $T_{k}(n)$ is also the unique graph with the maximum number of $K_t$ subgraphs among all $K_{k+1}$-free graphs.
\begin{thm}[Zykov~\cite{Zykov}]\label{erdosk}
Let $k$ and $t$ be integers such that $t < k+1$. Then for all $n$, the Tur\'{a}n graph $T_{k}(n)$ is the unique $K_{k+1}$-free graph on $n$ vertices containing the maximum number of $K_t$ subgraphs.
\end{thm}
We will need the following corollary in our calculations. 
\begin{cor}\label{zykov_K5_corollary}
Let $G$ be a $K_{k+1}$-free graph on $n$ vertices. Then 
\[ \nu(K_5,G) \leq \frac{k^4-10k^3+35k^2-50k+24}{k^4} n^5 + o(n^5).\]
\end{cor}

Alon and Shikhelman~\cite{Alon} proved the following analogue of the K\H{o}v\'{a}ri-S\'{o}s-Tur\'{a}n Theorem:
$$\ext(n,K_3,K_{s,t}) = O(n^{3-3/s}).$$
They also proved that for fixed integers $t < k$, if $F$ is a $k$-chromatic graph:
$$\ext(n,K_t,F) = \binom{k-1}{t}\left(\frac{n}{k-1}\right)^t + o(n^t).$$

In~\cite{GyoriPach1991}, Gy\H{o}ri, Pach, and Simonovits studied a handful of cases where $F = K_r$. The order of magnitude of $\ext(n,C_k,C_{\ell})$ is known for all $\ell \geq 3$ and $k \geq 3$, see Gishboliner and Shapira~\cite{GShapira}. 
The asymptotic value of $\ext(n,C_k,C_4)$ was determined by Gerbner, Gy\H{o}ri, Methuku, and Vizer~\cite{gerbner2017generalized}. They proved a variety of results on $\ext(n,F,H)$ when $F$ and $H$ were both cycles. 
This includes showing that $\ext(n,C_{2\ell},C_{2k}) = \Theta(n^{\ell})$ for $k,\ell \geq 2$ and \[\ext(n,C_4,C_{2k}) = (1+o(1))\frac{(k-1)(k-2)}{4}n^2\]
for $k \geq 2$.
In \cite{palmer}, Gerbner and Palmer provided more general bounds on $\ext(n,H,F)$. In particular, they showed that if $H$ and $F$ are graphs and $\chi(F) = k$, then 

$$ \ext(n,H,F) \leq \ext(n,H,K_k) + o\left(n^{|H|}\right).$$
Additionally, they extended the result of Gishboliner and Shapira to show that for all $k$ and $t$,
\[
\ext(n,C_k,K_{2,t}) = \left( \frac{1}{2k} +o(1) \right)(t-1)^{k/2}n^{k/2},
\]
and
\[
\ext(n,P_k,K_{2,t}) = \left( \frac{1}{2} +o(1) \right)(t-1)^{(k-1)/2}n^{(k+1)/2},
\]
where $P_k$ is a path on $k$ vertices.

In~\cite{stars}, Cutler, Nir, and Radcliffe determined the asymptotic value of $\ext(n,S_{t},K_{k+1})$, where $S_{t}$ is the star with $t$ leaves. In particular, they showed that while the extremal graph must be complete multi-partite, it is not always isomorphic to the Tur\'{a}n graph $T_k(n)$.  The study of the function $\ext(n,K_3,H)$ has seen recent attention as well. In particular, the function $\ext(n,K_3,C_5)$ was studied in \cite{Alon,bollobas2008,ergemlizdeC5free}. In \cite{mubayi2020triangles}, Mubayi and Mukherjee studied the function $\ext(n,K_3,H)$ for a handful of other $3$-chromatic graphs $H$. 

In~\cite{gerbner2020exact}, Gerbner and Palmer found a handful of cases where the value of $\ext(n,H,F)$ is achieved by the Tur\'{a}n graph and in~\cite{gerbner2020generalized}, Gerbner studied the function $\ext(n,H,F)$ when $H$ and $F$ each have at most $4$ vertices. 
Recently, the authors of~\cite{huynh2020subgraph} studied the problem of maximizing the number of copies of a graph $H$ in some graph $G$ embedded in a particular surface.

In 1984, Erd\H{o}s conjectured that the balanced blow-up of $C_5$ on $n$ vertices maximizes the number of five-cycles among all triangle-free graphs of order $n$. If $G$ is a graph on $m$ vertices, then the \emph{balanced blow-up} of $G$ on $n$ vertices is the graph $G(n)$ obtained from $G$ by replacing each vertex of $G$ with an independent set of size $\left\lfloor \frac{n}{m} \right\rfloor$ or $\left\lceil \frac{n}{m} \right\rceil$, and replacing each edge in $G$ with a complete bipartite graph on the corresponding sets. The problem of determining $\ext(n,C_5,K_3)$ was known as the Pentagon Problem of Erd\H{o}s. In a sense, a graph with $\ext(n,C_5,K_3)$ five-cycles is the ``least bipartite" triangle-free graph on $n$ vertices when measured by the number of 5-cycles. In posing this question, Erd\H{o}s also proposed the following two measures of ``non-bipartiteness''~\cite{Erdospent}.
\begin{enumerate}
    \item  The minimal possible number of edges in a subgraph spanned by half the vertices.
    
    \item The minimal possible number of edges that have to be removed to make the graph bipartite, which is equivalent to the problem of max cut.
\end{enumerate}
In 1989, Gy\H{o}ri~\cite{Gyoripent} showed that a triangle-free graph on $n$ vertices contains at most  $1.03\left(\frac{n}{5}\right)^5$ five-cycles. 
In 2012, Grzesik \cite{Grzesik} and independently in 2013, Hatami Hladk\'{y}, Kr\'{a}l', Norin, and Razborov \cite{Hatami} showed that a triangle-free graph on $n$ vertices contains at most $\left(\frac{n}{5}\right)^5 + o(n^5)$ five cycles. 
Moreover, a matching lower bound is given by the balanced blow-up of $C_5$ when $n$ is divisible by $5$. The authors of~\cite{Hatami} also proved that for large enough $n$, the balanced blow-up of a $C_5$ on $n$ vertices is the unique extremal graph. In 2018, Lidick\'{y} and Pfender \cite{Lidicky} proved that a balanced $C_5$ blow-up is the unique extremal construction for all $n$, with the exception of $n=8$. This observation was made by Michael~\cite{mobius}, who showed that the M\"{o}bius ladder on $8$ vertices contains the same number of five cycles as the balanced $C_5$ blow-up. 

Palmer~\cite{palmertalk} suggested a generalization to the Pentagon Problem of Erd\H{o}s: maximizing the number of five-cycles in $K_{k+1}$-free graphs for $k \geq 3.$ Observe that in the more general case, the problem of maximizing the number of non-induced $C_5$ subgraphs is different from maximizing the number of induced $C_5$ subgraphs.

In this paper, we will discuss the non-induced case. 
Let $H$ and $G$ be graphs on $n_1$ and $n_2$ vertices, respectively. The \emph{density} $d(H,G)$ of $H$ in $G$ is given by
$$d(H,G) = \nu(H,G)\binom{n_2}{n_1}^{-1}.$$
Normally, $n_2^{-n_1}$ would be used as the scaling factor for defining the density of $H$ in $G$. We will use $\binom{n_2}{n_1}^{-1}$, since this is more natural in proofs involving the flag algebra method. Let 
\begin{equation}\label{Opt definition}
\OPT_k(C_5) = \lim\limits_{n \to \infty} \max\limits_{G_n \in \mathcal{F}_n^k}d(C_5,G_n),
\end{equation}
where $\mathcal{F}_n^k$ is the set of all $K_{k+1}$-free graphs on $n$ vertices. Note that since $d(C_5,G)$ measures the density of non-induced $C_5$ subgraphs in a graph $G$, this parameter will often have a value greater than one.
For example, $d(C_5,K_\ell) = 12$ for all $\ell \geq 5$.
Our main goal is to prove the following theorem. 

\begin{thm}\label{main}
Let $k \geq 3$ be an integer. Then 
\begin{enumerate}
    \item[(i)] $\OPT_k(C_5) = \frac{1}{k^4}(12k^4 - 60k^3 + 120k^2 - 120k + 48).$
    \item[(ii)] If $n$ is sufficiently large,
    then $T_k(n)$ is the unique $K_{k+1}$-free graph on $n$ vertices for which $\nu(C_5,T_k(n)) = \ext(n,C_5,K_{k+1})$. 
\end{enumerate}
 \end{thm}

Since our result forbids $(k+1)$-cliques, Tur\'{a}n's Theorem implies that the number of edges in an extremal graph cannot be more than in $T_k(n)$. Interestingly, the authors of~\cite{BennetC5lowerbound} proved that if $G$ is a graph with at least $\frac{k-1}{k}\binom{n}{2}$ edges, then the Tur\'{a}n graph provides a lower bound on the number of five-cycles contained in $G$.

The proof of Theorem~\ref{main}(i) uses flag algebras to calculate the upper bound for $\OPT_k(C_5)$. The second part is done by stability and exact structure arguments. Unlike typical applications of the flag algebra method, our result does not need computer assistance for the calculations involving flag algebras.
However, it is still convenient to use a computer for the purpose of multiplying and expanding rational functions. 

In the next section, we will give a brief overview of the flag algebra method. Section 2 contains the proof of Theorem~\ref{main}(i). Then we prove a stability lemma in Section 3, and use it to prove Theorem~\ref{main}(ii) in Section 4. We will end with some concluding remarks and conjectures concerning the general behavior of the function $\ext(n,H,F)$. 

\subsection{The Flag Algebra Method}

Introduced by Razborov~\cite{Razborov}, the flag algebra method provides a framework for computationally solving problems in extremal combinatorics. Flag algebras have been used to solve problems on hypergraphs~\cite{hyperjump,hyperflag,hyperflag2,hyperflag3}, permutations~\cite{permutationflag}, graph decomposition problems~\cite{decompflag}, and oriented graphs~\cite{orientflag} among many other applications. 
Here we will give a brief introduction and description of the notation and theory we will need for our result.
We will not prove any claims since they have already been proven by Razborov~\cite{Razborov}. Another overview of flag algebras can be found in \cite{Flagfirst}. 

Let $H$ and $G$ be graphs on $n_1$ and $n_2$ vertices, respectively, such that $n_1 \leq n_2$. If $X \subseteq V(G)$, we will denote the induced subgraph of $G$ on the vertices of $X$ by $G[X]$. Let a subset $X$ be selected uniformly at random from $V(G)$ such that $|X| = n_1$. Then $P(H,G)$ is the probability that $G[X]$ is isomorphic to $H$. 

A sequence of graphs $(G_n)_{n \geq 1}$ of increasing order is said to be \emph{convergent} if for every finite graph $H$, the following limit converges:
    $$\lim\limits_{n \to \infty} P(H,G_n).$$
Let $\F$ denote the set of all graphs up to isomorphism, and let $\F_{\ell}$ denote the set of all graphs on $\ell$ vertices up to isomorphism. Let $\R \F$ denote the set of all formal finite linear combinations of graphs in $\F$. A \emph{type of size k} is a graph $\sigma$ on $k$ labelled vertices labeled by $[k]=\{1,\ldots,k\}$. If $\sigma$ is a type of size $k$ and $F$ is a graph on at least $k$ vertices, then an \emph{embedding} of $\sigma$ into $F$ is an injective function $\theta: [k] \to V(F)$, such that $\theta$ gives an isomorphism between $\sigma$ and $F[\im( \theta)]$. A $\sigma$-flag is a pair $(F,\theta)$ where $F$ is a graph and $\theta$ is an injective function from $[k]$ to $V(F)$ that defines a graph isomorphism of $F[\im(\theta)]$ and $\sigma$. In this way, $\sigma$ can be thought of as a labelled subgraph of $F$. Two $\sigma$-flags $F$ and $G$ are isomorphic if there exists a graph isomorphism between $F$ and $G$ that preserves the labelled subgraph $\sigma$. 

Let $\F^{\sigma}$ denote the set of all $\sigma$-flags and $\F_{\ell}^{\sigma}$ denote the set of all $\sigma$-flags on $\ell$ vertices. Observe that if $\sigma$ is the empty graph, then $\F^{\sigma} = \F$. For two $\sigma$-flags $F$ and $G$ with $|V(F)| \leq |V(G)|$, let $P(F,G)$ denote the probability that an injective map from $V(F)$ to $V(G)$ that fixes the labeled graph $\sigma$ induces a copy of $F$ in $G$. Razborov showed that there exists an algebra $\mathcal{A^{\sigma}}$ after some factorization of $\R\F^{\sigma}$. In doing so, he defined addition and multiplication on the elements of $\R\F^{\sigma}$. Addition can be defined in the natural way, by simply adding coefficients of the elements in $\R\F^{\sigma}$. We will now describe how to define multiplication of elements in $\mathcal{A}^{\sigma}$. 

Let $(G,\theta) \in \F^{\sigma}$ be a $\sigma$-flag on $n$ vertices. Let $(F_1,\theta_1),(F_2,\theta_2) \in \F^{\sigma}$ be two $\sigma$-flags for which $|V(F_1)| + |V(F_2)| \leq n + |V(\sigma)|$. Let $X_1$ and $X_2$ be two disjoint sets of sizes $|V(F_1)| - |V(\sigma)|$ and $|V(F_2)|-|V(\sigma)|$ respectively, selected uniformly at random from $V(G)\setminus\im(\theta)$. We will define the \emph{density of $F_1$ and $F_2$ in $G$}, denoted $P(F_1,F_2;G)$ as the probability that 
$(G[X_1 \cup \im(\theta)],\theta)$ is isomorphic to $(F_1,\theta_1)$ and $(G[X_2 \cup \im(\theta)],\theta)$ is isomorphic to $(F_2,\theta_2)$. 

It can be shown that as $n$ grows, then the density of $F_1$ and $F_2$ is approximately equal to the product of their individual densities: 
\begin{align}\label{eq:error}
    |P(F_1,F_2;G) - P(F_1,G)P(F_2,G)| \leq O(n^{-1}).
\end{align}
Given this fact, if $|V(F_1)| + |V(F_2)| = \ell - |V(\sigma)|$ we could ideally define multiplication in $\mathcal{A}^{\sigma}$ by 
\begin{align}\label{eq:multiplication}
F_1 \cdot F_2 = \sum\limits_{F \in \F_{\ell}^\sigma} P(F_1;F_2;F) F.
\end{align}
The issue with this, however, is that the product $F_1 \cdot F_2$ could be also written as a linear combination of elements in $\mathcal{F}^\sigma_{\ell'}$ for any $\ell'>\ell$. Hence, before defining $\mathcal{A}^{\sigma}$ we factor out all expressions of the form
\begin{equation}\label{flagfactor}
F - \sum\limits_{F' \in \F_{\ell}^\sigma} P(F,F')F'
\end{equation}
from $\R\F^{\sigma}$.
Note that \eqref{flagfactor} corresponds to the law of total probability and hence it should behave as 0 when added to another linear combination. 
Let $\mathcal{K}^{\sigma}$ be the linear subspace of $\R\F^{\sigma}$ containing all expressions of the form (\ref{flagfactor}). We define $\mathcal{A}^{\sigma}$ to be $\R\F^{\sigma}$ factorized by $\mathcal{K}^{\sigma}$, and we define multiplication in $\mathcal{A}^{\sigma}$ by naturally extending \eqref{eq:multiplication}.  

Returning to the idea of convergent sequences of graphs, let $\Hom^+(\mathcal{A}^{\sigma},\R)$ be the set of all homomorphisms from $\mathcal{A}^{\sigma}$ to $\R$ such that for each $\phi \in \Hom^+(\mathcal{A}^{\sigma},\R)$ and $H \in \mathcal{F}^{\sigma}$, $\phi(H) \geq 0$. 
If $\sigma$ has order 0, we omit it in the notation.
Razborov showed that each homomorphism in  $\Hom^+(\mathcal{A},\R)$ corresponds to some convergent graph sequence and vice versa~\cite{Razborov}. 

In any fixed graph $G$, we can express $d(C_5,G)$ as the sum of induced densities in the following way:

\begin{align}
d(C_5,G) = \sum\limits_{F_i \in \F_5}c_{F_i}^{C_5}P(F_i,G),\label{fullC5}
\end{align}
where $c_{F_i}^{C_5} = \nu(C_5,F_i)$. Hence, for any sequence of unlabelled graphs $(G_n)_{n \geq 1}$ and its corresponding homomorphism $\phi \in \Hom^+(\mathcal{A},\R)$, 

\begin{equation}\label{flaaaag}
    \lim\limits_{n \to \infty} d(C_5,G_n) = \sum\limits_{F_i \in \F_5}c_{F_i}^{C_5}\phi(F_i).
\end{equation}
Quite often in our computations to simplify notation, we will drop the function notation and simply write $F_i$ or draw the graph $F_i$ in place of $\phi(F_i)$. Under this notation equation~(\ref{flaaaag}) would be
\begin{equation*}
\lim\limits_{n \to \infty} d(C_5,G_n) = \sum\limits_{F_i \in \F_5}c_{F_i}^{C_5}F_i. 
\end{equation*}


Finally, while we will often work with $\sigma$-flags where $\sigma$ is not empty, flag algebras are often applied to questions concerning unlabelled graphs. In order to translate information from $\Hom^+(\mathcal{A}^{\sigma},\R)$ to $\Hom^+(\mathcal{A},\R)$
Razborov defined the \emph{unlabelling operator} which is a linear operator 
$$\llbracket \cdot \rrbracket_{\sigma}:\R\F^{\sigma} \to \R\F,$$
such that for any $\sigma$-flag $F=(H,\theta)$, $\llbracket F \rrbracket_{\sigma} = q_{\sigma}(F)H$, where $q_{\sigma}(F)$ is equal to the probability $(H,\theta')$ is isomorphic to $F$ where $\theta'$ is a randomly chosen injective mapping $\theta'$ from $[k]$ to $V(H)$. It can be shown that for any $a \in \mathcal{A}^{\sigma}$ and any $\phi \in \Hom^+(\mathcal{A},\R)$,
\begin{align}\label{eq:square}
\phi \left( \llbracket a \cdot a \rrbracket_{\sigma}\right) \geq 0.
\end{align}
We will frequently make use of this fact in our computations. 

If a flag algebra calculation has a constant number of terms and operations, then it can be interpreted as a calculation in a graph of order $n$ with an error term $O(n^{-1})$ coming from~\eqref{eq:error}.

\section{Proof of Theorem~\ref{main}(i)}
In this section we will prove Theorem~\ref{main}(i). First we will provide a lower bound by counting the number of five cycles in the Tur\'{a}n graph. Next, using the flag algebra method, we will provide a matching upper bound. The proof of the upper bound when $k = 3$ is slightly different than the proof when $k \geq 4$.  

\medskip\noindent 
\begin{proof}[Proof of Theorem~\ref{main}(i)]
First we will count the number of five-cycles in $T_k(n)$, which will give an asymptotic lower bound.  
Observe that the only induced subgraphs of $T_k(n)$ on five vertices containing a five-cycle are $\Fhhsml, \Fggsml$, and $\Fffsml$. 
There are $\binom{k}{5}\left(\frac{n}{k}\right)^5$ copies of $\Fhhsml$ in $T_k(n)$, with every such graph containing $12$ distinct $C_5$ subgraphs. 
There are $4\binom{k}{4}\binom{n/k}{2}\left(\frac{n}{k}\right)^3$ copies of $\Fggsml$ in $T_k(n)$, with every such graph containing $6$ distinct $C_5$ subgraphs. 
Finally, there are $3\binom{k}{3}\binom{n/k}{2}^2\frac{n}{k}$ copies of $\Fffsml$ in $T_k(n)$, with every such graph containing $4$ distinct $C_5$ subgraphs.
Thus, 
\[
\nu(C_5,T_k(n)) = 12\binom{k}{5}\left(\frac{n}{k}\right)^5 + 24\binom{k}{4}\binom{n/k}{2}\left(\frac{n}{k}\right)^3 + 12\binom{k}{3}\binom{n/k}{2}^2\frac{n}{k} + o(n^5),
\]
where the error term $o(n^5)$ accounts for the cases where $n$ is not divisible by $k$. 
This implies that for all $n$, 
$$d(C_5,T_k(n)) \geq \nu(C_5,T_k(n))\binom{n}{5}^{-1} = \frac{1}{k^4}(12k^4 - 60k^3 + 120k^2 - 120k + 48) + o(1).$$

Now we will calculate an asymptotic upper bound. Unless it is stated otherwise, assume that $k \geq 3$. Let $\F_5 = \{F_0, \ldots ,F_{33}=K_5\}$ be the set of unlabeled graphs (up to isomorphism) on five vertices. Each of these graphs is pictured in Table~\ref{tab:fivevertexgrahs} in the Appendix.  
After removing each $c_{F_i}^{C_5}$ for which $c_{F_i}^{C_5} = 0$ from \eqref{fullC5} we get
\begin{equation}\label{easyC5}
\lim\limits_{n \to \infty} d(C_5,G_n) = \Faa + \Fbb + \cdot \Fcc + 2\cdot \Fdd + 2\cdot \Fee + 4\cdot \Fff + 6\cdot \Fgg + 12\cdot \Fhh. 
\end{equation}
Since $\mathcal{F}_5$ contains all graphs on five vertices (up to isomorphism) and $\sum_{i=0}^{33} F_i = 1$,
$$\lim_{n \to \infty} d(C_5,G_n) \leq \max \{c_{F_i}^{C_5} : F_i \in \F_5\}.$$ 
Therefore, 
$$\OPT_k(C_5) \leq \max\{c_{F_i}^{C_5} : F_i \in \F_5\}.$$
Given this fact, our goal is to find appropriate constants $c_{F_i}$ so that
$$\lim\limits_{n \to \infty} d(C_5,G_n) \leq \sum\limits_{F_i \in \F_5}c_{F_i}F_i,$$
and $\max \{c_{F_i} : F_i \in \F_5\}$ is as small as possible. To do so, we can take advantage of properties that we know must be true of all $\phi \in \Hom^+(\mathcal{A},\R)$ which correspond to $K_{k+1}$-free convergent sequences of graphs. Additionally, using labeled flags, we can derive nonnegative expressions of unlabeled graphs in $\mathcal{F}_5$. We define
\begin{align*}
\sigma_1 &= \threesigone & \sigma_2 &= \threesigtwo & \sigma_3 &= \threesigthree
\end{align*}
so that $\F_4^{\sigma_1}$, $\F_4^{\sigma_2}$, and $\F_4^{\sigma_3}$ denote three sets of labeled flags on four vertices. 
By 
\eqref{eq:square}, the following expressions are nonnegative for all $k \geq 3$. 

\begin{enumerate}
    \setlength\itemsep{2 em}
    \item $P_1(k) = 10 \cdot \left\llbracket \left((k-1)\fourone -\fourtwo\right)^2 \right\rrbracket_{\sigma_1}  = $
    
    $(10k^2-20k+10)\Fa + (k^2-2k+1)\Fb + (-k+1)\Fd + (-4k + 4)\Fe + \Fs + \Ft$
    
    \item $P_2(k) = 30 \cdot \left\llbracket \left( (k-2)\fourthree - \fourfour \right)^2 \right\rrbracket_{\sigma_2} = $
    
    $(3k^2-12k+12)\Fs + (k^2-6k+8)\Fdd + (-4k+10)\Fff + 3\Fgg$

   \item  $P_3(k) = 30 \cdot \left\llbracket \left( (k-2)\fourfive - \fourfour \right)^2 \right\rrbracket_{\sigma_2} = $
   
    $(6k^2 - 24k + 24)\Fe + (k^2-4k+4)\Fl + (-k+2)\Fr + (-6k+12)\Ft + 2\Fff + 3\Fgg$
    
   \item  $P_4(k) = 30 \cdot \left\llbracket \left( \foursix - \fourseven \right)^2 \right\rrbracket_{\sigma_3} = $
   
       $ 6\Ft - \Fcc + 2\Fee -4\Fff$
        
    \item $P_5(k) = 30 \cdot \left\llbracket \left( (k-3)\foursix + (k-3)\fourseven - 2\foureight \right)^2 \right\rrbracket_{\sigma_3} = $

    $(6k^2-36k+54)\Ft + (2k^2-20k+42)\Fee + (4k^2 - 24k + 36)\Fff + (-24k + 84)\Fgg + 120\Fhh$
\end{enumerate}
Additionally, we can apply Corollary~\ref{zykov_K5_corollary}, which states that for any $K_{k+1}$-free convergent sequence of graphs, 

$$\Fhh \leq \frac{k^4-10k^3+35k^2-50k+24}{k^4}.$$

At this point we will split the proof into the two cases where $k \geq 4$ and $k = 3$. To gain some intuition as to why this is necessary, we can consider the previous inequality. When $k = 3$ or $k = 4$, the previous bound implies that $\Fhhsml = 0$. The issue is that it does not give any information about the density of $K_4$, which is also equal to zero when $k = 3$. Thus, two slightly different proofs are required for $k = 3$ and $k \geq 4$. 

\medskip
\textbf{Case 1:} Suppose that $k \geq 4.$ Since $\sum_{i=0}^{33} F_i = 1$,

\begin{equation}\label{K5bound}
\Fhh \leq \sum\limits_{i = 0}^{33}F_i\left(\frac{k^4-10k^3+35k^2-50k+24}{k^4}\right).
\end{equation}
After rearranging the terms from (\ref{K5bound}) we obtain the following constraint on the elements of $\F_5$:
\begin{equation}\label{K5boundsecond}
0 \leq \sum\limits_{i = 0}^{32}F_i \cdot \frac{k^4-10k^3+35k^2-50k+24}{k^4} + \Fhh \cdot \frac{-10k^3 + 35k^2 -50k+24}{k^4}.
\end{equation}
Let
$$Z(k) =\sum\limits_{i = 0}^{32}F_i \cdot \frac{k^4-10k^3+35k^2-50k+24}{k^4} + \Fhh \cdot \frac{-10k^3 + 35k^2 -50k+24}{k^4}.$$

It is straightforward to verify that the following rational functions are nonnegative for all $k \geq 4$. 
\begin{align*}
\textstyle
 p_1(k) &= 
 \textstyle
 \frac{3(k^5 - 8k^4 + 22k^3 - 24k^2 + 8k)}{5k^7 - 35k^6 + 75k^5 - 48k^4}
 &
 p_2(k) &= 
 \textstyle
 \frac{10k^5 - 60k^4 + 109k^3 - 76k^2 + 18k}{5k^7 - 35k^6 + 75k^5 - 48k^4}
 \\
 p_3(k) &= 
 \textstyle
 \frac{5k^5 - 28k^4 + 45k^3 - 28k^2 + 6k}{5k^7 - 35k^6 + 75k^5 - 48k^4}
 &
p_4(k) &= 
\textstyle
\frac{5k^7 - 30k^6 + 53k^5 - 52k^4 + 94k^3 - 96k^2 + 24k}{4(5k^7 - 35k^6 + 75k^5 - 48k^4)}
\\
p_5(k) &= 
\textstyle
\frac{15k^5 - 60k^4 + 78k^3 - 40k^2 + 8k}{4(5k^7 - 35k^6 + 75k^5 - 48k^4)}
&
z(k) &= 
\textstyle
\frac{6(5k^3 - 20k^2 + 30k - 16)}{5k^3 - 35k^2 + 75k - 48}
\end{align*}
Thus, for any convergent sequence $(G_n)_{n \geq 1}$ of $K_{k+1}$-free graphs with $k \geq 4$, 
\begin{align}\label{C5final}
\lim\limits_{n \to \infty} d(C_5,G_n) \leq \sum\limits_{i = 0}^{33} c_{F_i}^{C_5}F_i + z(k)Z(k) + \sum\limits_{j=1}^5 p_j(k)P_j(k) = \sum\limits_{i = 0}^{33} c_{F_i}F_i,
\end{align}
where $c_{F_i}$ is the coefficient of $F_i$ after all expansions. 
Then 
$$\OPT_k(C_5) \leq \max\{c_{F_i}: F_i \in \F_5\}.$$
The values of $c_{F_i}$ for each $F_i \in \F_5$ are listed below. 

\begin{itemize}
 \item $C_1(k) = c_{\Fasml} = c_{\Fesml} = c_{\Fssml} = c_{\Ftsml} $
 
 \medskip
 $= c_{\Fffsml} = c_{\Fggsml} = c_{\Fhhsml} = \frac{60k^7 - 720k^6 + 3600k^5 - 9876k^4 + 16320k^3 - 16440k^2 + 9360k - 2304}{5k^7 - 35k^6 + 75k^5 - 48k^4}.$

\item $C_2(k) = c_{\Fbsml} =\frac{33k^7-450k^6+2547k^5-7824k^4+14214k^3-15360k^2+9144k-2304}{5k^7 - 35k^6 + 75k^5 - 48k^4}.$

\item $C_3(k) = c_{\Fcsml} = c_{\Ffsml} = c_{\Fgsml} = c_{\Fhsml} = c_{\Fisml} = c_{\Fjsml} = c_{\Fksml} = c_{\Fmsml} = c_{\Fnsml} = c_{\Fosml} = c_{\Fpsml} = c_{\Fqsml} = c_{\Fusml}$

\medskip
$ = c_{\Fvsml} = c_{\Fwsml} = c_{\Fxsml} = c_{\Fysml} = c_{\Fzsml} = \frac{30k^7 - 420k^6 + 2430k^5 - 7596k^4 + 13980k^3 - 15240k^2 + 9120k - 2304}{5k^7 - 35k^6 + 75k^5 - 48k^4}.$

\item $C_4(k) = c_{\Fdsml} = \frac{30k^7 - 423k^6 + 2457k^5 - 7686k^4 + 14118k^3 - 15336k^2 + 9144k - 2304}{5k^7 - 35k^6 + 75k^5 - 48k^4}.$

\item $C_5(k) = c_{\Flsml} = \frac{35k^7 - 468k^6 + 2607k^5 - 7916k^4 + 14278k^3 - 15367k^2 + 9144k - 2304}{5k^7 - 35k^6 + 75k^5 - 48k^4}.$

\item $C_6(k) = c_{\Frsml} = \frac{30k^7 - 425k^6 + 2468k^5 - 7697k^4 + 14098k^3 - 15302k^2 + 9132k - 2304}{5k^7 - 35k^6 + 75k^5 - 48k^4}.$

\item $C_7(k) = c_{\Faasml} = c_{\Fbbsml}= \frac{35k^7 - 455k^6 + 2505k^5 - 7644k^4 + 13980k^3 - 15240k^2 + 9120k - 2304}{5k^7 - 35k^6 + 75k^5 - 48k^4}.$

\item $C_8(k) = c_{\Fccsml} = \frac{(135/4)k^7 - (895/2)k^6 + (9967/4)k^5 - 7631k^4 + (27913/2)k^3 - 15216k^2 + 9114k - 2304}{5k^7 - 35k^6 + 75k^5 - 48k^4}.$

\item $C_{9}(k) = c_{\Fddsml} = \frac{50k^7 - 610k^6 + 3129k^5 - 8902k^4 + 15326k^3 - 15956k^2 + 9264k - 2304}{5k^7 - 35k^6 + 75k^5 - 48k^4}.$

\item $C_{10}(k) = c_{\Feesml} = \frac{50k^7 - 610k^6 + 3103k^5 - 8758k^4 + 15050k^3 - 15748k^2 + 9216k - 2304}{5k^7 - 35k^6 + 75k^5 - 48k^4}.$
\end{itemize}
\begin{claim}\label{claim that C1 is max}
For all $i = 1,\dots,10$ and $k \geq 4$, $C_1(k) \geq C_i(k)$. 
\end{claim}
\begin{proof}
Observe that for all $i = 1,\dots,10$, each polynomial $C_i(k)$ has the same denominator of $5k^7 -35k^5 + 75k^4 -48k^3$. It is straightforward to verify that $5k^7 -35k^5 + 75k^4 -48k^3$ is positive for all $k \geq 4$. By examining the leading coefficients in the numerator of each polynomial, it is straightforward to check that $C_1(k)$ is the largest for $k > 1000$. For $4 \leq k \leq 1000$, we have provided Sage code used to verify the claim in Appendix~\ref{proof that C1 is max}. 
\end{proof}

\medskip
By factoring $C_1(k)$ it follows that
\[
\OPT_k(C_5) \leq C_1(k) = \frac{1}{k^4}(12k^4 - 60k^3 + 120k^2 - 120k + 48),
\]
completing the proof of Theorem~\ref{main}(i) when $k \geq 4$.

\medskip
\textbf{Case 2:} Suppose that $k = 3$. Assume that $(G_n)_{n \geq 1}$ is a $K_4$-free convergent sequence of graphs. Each graph in the set $\mathcal{H}$ given below has a limit density equal to zero, and therefore can be removed from our calculations. 
    \begin{figure}[H]
    \centering
    $\mathcal{H} = \left\{
     \Fu,
     \hspace{1 cm}
     \Fy,
     \hspace{1 cm}
     \Fee,
     \hspace{1 cm}
     \Fgg,
     \hspace{1 cm}
     \Fhh
     \right\}$.
    \label{forbidden when k is three}
    \end{figure}
In this case, we will use the same polynomials $P_i(k)$ for $i = 1,2,3,4$ that were provided earlier in the proof. We will define one new polynomial $P_6$, which is nonnegative by \eqref{eq:square}. 

\[P_6 = \left\llbracket \left( \fourforkthree + \fourseven - \diamondforkthree \right)^2 \right\rrbracket_{\sigma_3} = \]
\[\Fl + 2 \cdot \Fz  - \Fx - 2 \cdot \Fdd + \Fr + 6 \cdot \Ft - 4 \cdot \Fff \geq 0\]
Now suppose that 
\begin{align*}
    p_1 &= 1/27   &
    p_2 &= 13/27  &
    p_3 &= 8/27   \\
    p_4 &= 2/9    &
    p_6 &= 17/54.  
\end{align*}
Then
$$d(C_5) = \lim\limits_{n \to \infty} d(C_5,G_n) \leq \sum\limits_{F \in \mathcal{F}_5 \setminus \mathcal{H}} \nu(C_5,F) F + \sum\limits_{j = 1}^4 p_jP_j(3) + p_6P_6.$$
Let $c_F$ denote the coefficient of each graph $F$ after combining each of the two sums. We provide the values of $c_F$ for each graph in $\mathcal{F}_5 \setminus \mathcal{H}$ for which $c_F \neq 0$. 
\begin{itemize}
    \item $c_{\Fasml} = c_{\Fesml} = c_{\Fssml} = c_{\Ftsml} = c_{\Fffsml} = \frac{40}{27}$
    \item $c_{\Fbsml} = \frac{4}{27}$
    \item $c_{\Fdsml} = -\frac{2}{27}$
    \item $c_{\Flsml} = \frac{11}{18}$
    \item $c_{\Frsml} = \frac{1}{54}$
    \item $c_{\Fxsml} = -\frac{17}{54}$
    \item $c_{\Faasml} = c_{\Fbbsml} = 1$
    \item $c_{\Fccsml} = \frac{7}{9}$
    \item $c_{\Fddsml} = \frac{8}{9}$
\end{itemize}
SageMath code to verify this calculation can be found in Appendix~\ref{proof that C1 is max}. It is straightforward to verify that
$$ d(C_5) \leq \max\{c_{F} : F \in \mathcal{F}_5 \setminus \mathcal{H}\} = \frac{40}{27}.$$
Furthermore, the set $T_3$ of graphs for which $c_{F_i} = \frac{40}{27}$ is given below.
    \begin{align}
    T_3 = \left\{
     \Fa,
     \hspace{1 cm}
     \Fe,
     \hspace{1 cm}
     \Fs,
     \hspace{1 cm}
     \Ft,
     \hspace{1 cm}
     \Fff
     \right\}
     \label{T3}
     \end{align}
\noindent This completes the proof of Theorem~\ref{main} (i).
\end{proof}

\subsection{Finding the Optimal Bound}
We will now give a short description on how we found the functions $z(k)Z(k)$ and $p_i(k)P_i(k)$ that were used in the proof of Theorem~\ref{main}(i).
If $F_j \in \mathcal{F}_5$ is a graph for which $c_{F_j} = \OPT_k(C_5)$, then we call $F_j$ a \emph{tight} subgraph. In our proof of Theorem~\ref{main}(i), the set $T$ given below contains the tight subgraphs when $k \geq 4$.
    \begin{figure}[H]
    \centering
    $T = \left\{
     \Fa,
     \hspace{1 cm}
     \Fe,
     \hspace{1 cm}
     \Fs,
     \hspace{1 cm}
     \Ft,
     \hspace{1 cm}
     \Fff,
     \hspace{1 cm}
     \Fgg,
     \hspace{1 cm}
     \Fhh
     \right\}$
    \end{figure}
The set $T_3$, defined in \eqref{T3}, contains the tight subgraphs when $k = 3$.
Note that $T_3$ are the $K_4$-free graphs from $T$.
The following lemma, which appears as Lemma 2.4.3 in~\cite{baberthesis}, states that any graph appearing with positive probability in the limit of $(G_n)_{n \geq 1}$ must be tight if $G_n$ are extremal graphs. 

\begin{lemma}[\cite{baberthesis}]\label{tightsubs}
Given $(G_n)_{n \geq 1}$ a convergent sequence of $K_{k+1}$-free graphs of increasing order, such that $d(C_5,G_n) \to \OPT_k(C_5)$. Let $d(H,G_{\infty})$ be the value of $\lim\limits_{n \to \infty} d(H,G_n)$. Then for any exact solution, $ d(H,G_{\infty}) > 0$ implies that $H$ must be a tight subgraph.
\end{lemma}

Using semidefinite programming, we verified that the conjectured upper bound of $\OPT_k(C_5)$ was correct for small values of $k$. 
In doing so, we were able to guess the correct types and labelled flags to use. 
It was a greedy process and there may be simpler solutions.
This corresponds to the polynomials $P_i$ for $i = 1,\dots,6$. 
Note that each labelled flag is a four-vertex graph appearing in the Tur\'{a}n graph.  
Next, Lemma~\ref{tightsubs} implies that each $F_i \in \mathcal{F}_5$ that is a subgraph of the Tur\'{a}n graph must have the property that $c_{F_i} = \OPT_k(C_5)$. 
Given this fact, we used SageMath to solve for the correct polynomials $p_i(k)$ and $z(k)$. 
These agreed with the values calculated by the semidefinite program for small $k$.

\section{Stability}
In this section we will prove a stability lemma which states that for any $K_{k+1}$-free graph $G$ on a sufficient number of vertices, if $G$ contains ``close" to the extremal number of five-cycles, then $G$ can be made isomorphic to $T_k(n)$ by adding or deleting a small number of edges. 

\begin{prop}\label{bounds}
For two positive integers $x_1$ and $x_2$, if $x_1 \geq x_2 +2$, then 
\begin{enumerate}
\item $x_1x_2 < (x_1-1)(x_2+1)$
\item $x_1\binom{x_2}{2} < (x_2+1)\binom{x_1 - 1}{2}.$
\end{enumerate}
\end{prop}
\begin{proof}
The first inequality is clear from the equation below:
$$(x_1-1)(x_2+1) = x_1x_2 + (x_1-x_2) - 1 \geq x_1x_2 + 1.$$
Since $x_2 - 1 < x_1 - 2$, the second inequality follows immediately from the first inequality. 
\end{proof}

\medskip
The next proposition follows immediately from Lemma~\ref{move a vertex}, which we will prove next. 
\begin{prop}\label{setsize}
For any complete $k$-partite graph $H$ on $n$ vertices, $\nu(C_5,H)$ is maximized when the sizes of the partite sets are as equal as possible. 
\end{prop}

The following lemma will show that if $H$ is a complete $k$-partite graph with unbalanced partite sets, then we can always increase the number of five-cycles in $H$ by moving the vertices as to make $H$ more balanced. Throughout the proof, we will assume for each $i = 1,\dots,k$ that $|X_i| = x_i$. For a graph $G$ and a vertex $v \in V(G)$, let $\nu_{G}(v,C_5)$ denote the number of five-cycles in $G$ that contain $v$. 

\begin{lemma}\label{move a vertex}
Let $H$ be a complete $k$-partite graph with partite sets $X_1,\dots,X_k$. Suppose that for two integers $i$ and $j$, 
$$x_i \geq x_j + 2.$$
Let $H'$ be the graph obtained from $H$ by deleting a vertex in $X_i$ and duplicating a vertex in $X_j$. Then 
$$\nu(C_5,H') > \nu(C_5,H).$$
\end{lemma}
\begin{proof}
By symmetry we may assume that $i = 1$ and $j = 2$. 
Let $H'$ be obtained from $H$ by removing some vertex $v \in X_1$ and adding a new vertex $v'$ to $X_2$, where $v'$ is a duplicate of some vertex in $X_2$.
Letting $X_1',\dots,X_k'$ denote the new partite sets in $H'$, $|X_1'| = x_1 - 1$, $|X_2'| = x_2 + 1$, and $|X_q'| = x_q$ for each $q \in \{3,\dots,k\}$. 

Since the only five-cycles that have been deleted from $H$ are those containing $v$, we only need to show that $\nu_{H'}(v',C_5) > \nu_H(v,C_5)$. Additionally, there is a one-to-one correspondence between the five cycles in $H$ containing $v$ and no other vertices in $X_1\cup X_2$ and the five cycles in $H'$ containing $v'$ and no other vertices in $X_1' \cup X_2'$. Because of this, we can focus only on those five-cycles which contained $v$ and at least one other vertex in $X_1 \cup X_2$. 

Let $c(v,n_1,n_2)$ denote the number of five cycles in $H$ containing $v$ along with $n_1$ and $n_2$ vertices in $X_1$ and $X_2$, respectively. We define $c'(v',n_1,n_2)$ in an identical manner, but pertaining to $v'$ and $H'$. In order to show that $\nu(C_5,H') > \nu(C_5,H)$, it suffices to show the following,
\begin{enumerate}
    \item $c'(v',1,0) + c'(v',2,0) + c'(v',0,1) > c(v,0,1) + c(v,0,2) + c(v,1,0)$, and
    \item $c'(v',1,1) + c'(v',2,1) > c(v,1,1) + c(v,1,2).$
\end{enumerate}
We will prove each of these inequalities as two separate claims. Throughout the proof we will assume that $I = \{3,\dots,k\}$.
\begin{claim}\label{single vertex increase prop}
$c'(v',1,0) + c'(v',2,0) + c'(v',0,1) > c(v,0,1) + c(v,0,2) + c(v,1,0)$.
\end{claim}
\begin{proof}
Since $H$ is a complete $k$-partite graph,
\begin{align*}
&c(v,0,1) = 6x_2 \cdot \sum\limits_{\{i,j\} \in \binom{I}{2}} \left[\binom{x_i}{2}x_j + \binom{x_j}{2}x_i\right] + 12x_2\cdot \sum\limits_{\{i,j,h\} \in \binom{I}{3}} x_ix_jx_h,\\
&c(v,0,2) = 4\binom{x_2}{2} \cdot \sum\limits_{i \in I} \binom{x_i}{2} + 6\binom{x_2}{2}\cdot \sum\limits_{\{i,j\} \in \binom{I}{2}} x_ix_j, \text{ and } \\
&c(v,1,0) = 4(x_1-1) \cdot \sum\limits_{\{i,j\} \in \binom{I}{2}} \left[\binom{x_i}{2}x_j + \binom{x_j}{2}x_i\right] + 6x_1\cdot \sum\limits_{\{i,j,h\} \in \binom{I}{3}} x_ix_jx_h.
\end{align*}
By counting in similar way in $H'$, 
\begin{align*}
&c'(v',1,0) = 6(x_1 - 1) \cdot \sum\limits_{\{i,j\} \in \binom{I}{2}} \left[\binom{x_i}{2}x_j + \binom{x_j}{2}x_i\right] + 12(x_1 - 1)\cdot \sum\limits_{\{i,j,h\} \in \binom{I}{3}} x_ix_jx_h, \\
&c'(v',2,0) = 4\binom{x_1 - 1}{2} \cdot \sum\limits_{i \in I} \binom{x_i}{2} + 6\binom{x_1 -1}{2}\cdot \sum\limits_{\{i,j\} \in \binom{I}{2}} x_ix_j, \text{ and } \\
&c'(v',0,1) = 4x_2\cdot \sum\limits_{\{i,j\} \in I} \left[\binom{x_i}{2}x_j + \binom{x_j}{2}x_i\right] + 6x_2\cdot \sum\limits_{\{i,j,h\} \in I} x_ix_jx_h.
\end{align*}

Since $x_1 \geq x_2 +2$, it follows that $c'(v',2,0) > c'(v',0,2)$. Thus, it suffices to show that 
\[c(v,0,1) + c(v,1,0) \leq c'(v',0,1) + c'(v',1,0).\]
It is straightforward to verify that 
\small{
\[6\cdot \sum\limits_{\{i,j\} \in \binom{I}{2}} \left[\binom{x_i}{2}x_j + \binom{x_j}{2}x_i\right] + 12\cdot \sum\limits_{\{i,j,h\} \in \binom{I}{3}} x_ix_jx_h \geq 4\cdot \sum\limits_{\{i,j\} \in \binom{I}{2}} \left[\binom{x_i}{2}x_j + \binom{x_j}{2}x_i\right] + 6\cdot \sum\limits_{\{i,j,h\} \in \binom{I}{3}} x_ix_jx_h.\]}
This immediately implies that 
\[ c(v,0,1) - c(v,0,1) \leq c'(v',1,0) - c'(v',1,0),\]
which proves the claim.
\end{proof}
\begin{claim}\label{some in x1 and some in x2}
$c'(v',1,1) + c'(v',2,1) > c(v,1,1) + c(v,1,2).$
\end{claim}
\begin{proof}
For convenience, we will count $c(v,1,1) + c(v,1,2)$ in the following way:
\begin{equation}\label{claim 35 less}
c(v,1,1) + c(v,1,2) = x_1x_2f_{11} + \binom{x_2}{2}x_1f_{21},
\end{equation}
where $f_{pq}$ is a function independent of the values $x_1$ and $x_2$ used to count the number of five cycles containing $v$, $p$ vertices from $X_1$, and $q$ vertices from $X_2$. Using the same method to count $c'(v',1,1) + c'(v',2,1)$, we get
\small{
\begin{equation}\label{claim 35 more}
c'(v',1,1) + c'(v',2,1) = (x_1-1)(x_2+1)f_{11} + \binom{x_1-1}{2}(x_2+1)f_{12}. 
\end{equation}
}
By Proposition~\ref{bounds}, 
\[ (x_1-1)(x_2+1)f_{11} > x_1x_2f_{11}.\]
Moreover, since the sizes of each set $X_j$ for all $j \in I$ have not changed, $f_{12} = f_{21}$. Therefore, 
\[\binom{x_1-1}{2}(x_2+1)f_{12} > \binom{x_2}{2}x_1f_{21}\]
by Proposition~\ref{bounds}, completing the proof of the claim. \end{proof}

As each of Claims~\ref{single vertex increase prop} and \ref{some in x1 and some in x2} are true, it follows that  $\nu(C_5,H') > \nu(C_5,H)$, completing the proof of Lemma~\ref{move a vertex}. 
\end{proof}

\medskip
For two graphs $G$ and $H$ of the same order, let $\Dist(G,H)$ equal the minimum number of adjacencies that one needs to change in $G$ in order to obtain a graph isomorphic to $H$. The parameter $\Dist(G,H)$ is commonly known as the \emph{edit distance} between $G$ and $H$. Our main goal of this section is to prove the following lemma. 

\begin{lemma}[Stability Lemma]\label{friststab}
For every $\varepsilon > 0$, there exists an $n_0$ and $\varepsilon_F >0$ such that for every $K_{k+1}$-free graph $G$ of order $n \geq n_0$ with $d(C_5,G) \geq OPT_k(C_5) - \varepsilon_F,$ the edit distance between $G$ and $T_k(n)$ is at most $\varepsilon n^2$. 
\end{lemma}

The proof of Lemma~\ref{friststab} requires the following two lemmas along with Lemma~\ref{tightsubs}. For a family of graphs $\mathcal{F}$, we say that a graph $G$ is \emph{$\mathcal{F}$-free} if $G$ does not contain any member of $\mathcal{F}$ as an induced subgraph. 

\begin{lemma}[Induced Removal Lemma~\cite{removref}]\label{removlem}
Let $\mathcal{F}$ be a set of graphs. For each $\varepsilon > 0 $, there exist $n_0 \geq 0$ and $\delta > 0$ such that for every graph $G$ of order $n_0 \geq n$, if $G$ contains at most $\delta n^{|V(H)|}$ induced copies of $H$ for every $H \in \mathcal{F}$, then $G$ can be made $\mathcal{F}$-free by removing or adding at most $\varepsilon n^2$ edges from $G$. 
\end{lemma}

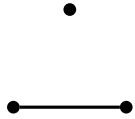
\begin{figure}[h]
\centering
\begin{tikzpicture}[scale = 0.75]
\draw node[inner sep=1.7pt, outer sep=0pt, circle, fill] (a) at (0,0) {};
\draw node[inner sep=1.7pt, outer sep=0pt, circle, fill] (b) at (-1,-1.732) {};
\draw node[inner sep=1.7pt, outer sep=0pt, circle, fill] (c) at (1,-1.732) {};
\draw[color=black,line width=1.2pt] (b)--(c);
\end{tikzpicture} 
\caption{$\overline{P_3}$}\label{fig:p3}
\end{figure}

Let $(G_n)_{n \geq 1}$ be a convergent sequence of $K_{k+1}$-free graphs. In the proof of Theorem~\ref{main}(i), we found constants $c_{F_i}$ for each $F_i \in \mathcal{F}_5$ such that
$$ d(C_5,G_n) \leq \sum\limits_{i = 0}^{33} c_{F_i} F_i \leq \max \{c_{F_i} : F_i \in \mathcal{F}_5\}$$
and
$$\max \{c_{F_i} : F_i \in \mathcal{F}_5\} = \OPT_k(C_5) = \frac{1}{k^4}(12k^4 - 60k^3 + 120k^2 - 120k + 48).$$
Let $\overline{P_3}$ be the three vertex graph with exactly one edge; see Figure~\ref{fig:p3}. The goal of Lemma~\ref{secstab} is to prove that if $\lim_{n \to_\infty}(C_5,G_n) = \OPT_k(C_5)$,
then $\lim_{n \to_\infty}(\overline{P_3},G_n) = 0$. 

\begin{lemma}\label{secstab}
For each $\delta_F > 0$, there exists $\varepsilon_F > 0$ and $n_0 = n_0(\delta_F)$ such that any $K_{k+1}$-free graph $G$ on $n \geq n_0$ vertices with $d(C_5,G) > OPT_k(C_5) - \varepsilon_F$ contains at most $\delta_F n^3$ induced copies of $\overline{P_3}$. 
\end{lemma}
\begin{proof}
Let $(G_n)_{n \geq 1}$ be a convergent sequence of $K_{k+1}$-free graphs maximizing the number of five-cycles. Let $T$ be the set of tight subgraphs in $\mathcal{F}_5$ given by the proof of Theorem~\ref{main}(i). This is the same set $T$ provided at the end of Section 2. 

    \begin{figure}[H]
    \centering
    $T = \left\{
     \Fa,
     \hspace{1 cm}
     \Fe,
     \hspace{1 cm}
     \Fs,
     \hspace{1 cm}
     \Ft,
     \hspace{1 cm}
     \Fff,
     \hspace{1 cm}
     \Fgg,
     \hspace{1 cm}
     \Fhh
     \right\}$.
    \end{figure}

Observe that for each graph $F \in T$, 
\[
P(\overline{P_3},F) = 0.
\]
Since $T$ contains the set of tight graphs, the following is a consequence of \eqref{flagfactor} for the sequence $(G_n)_{n\geq 1}$,
\[
\overline{P_3} = \sum\limits_{i = 0}^{33} P(\overline{P_3},F_i) F_i = 0.
\]
It follows that for the sequence $(G_n)_{n\geq 1}$,
\[
\lim\limits_{n \to \infty} d(\overline{P_3},G_n) = 0,
\]
which completes the proof of Lemma~\ref{secstab}. 
\end{proof}

\medskip
\begin{proof}[Proof of Lemma~\ref{friststab}]
Let $\varepsilon_I > 0$ and $\varepsilon_F > 0$, which we will determine later. 
By Lemma~\ref{removlem}, there exists a $\delta_F > 0$ and an $n_0$ such that any $K_{k+1}$-free graph on $n \geq n_0$ vertices containing at most $\delta_F n^3$ copies of $\overline{P_3}$ can be made $\overline{P_3}$-free after changing at most $\varepsilon_I n^2$ adjacencies. Assume that $G$ is a graph on $n \geq n_0$ vertices such that 
\[
d(C_5,G) > OPT_k(C_5) - \varepsilon_F,
\]
where $n_0$ is large enough to satisfy the conditions of Lemmas~\ref{removlem} and \ref{secstab} so that $G$ contains at most $\delta_Fn^3$ copies of $\overline{P_3}$. 
Moreover, for sufficiently small $\varepsilon_I$, 
\[
d(C_5,G) > OPT_{k-1} (C_5) + 2\cdot 5! \cdot \varepsilon_I.
\]
Using Lemma~\ref{removlem}, let 
$G'$ be a $\overline{P_3}$-free graph obtained from
$G$ by changing at most $\varepsilon_I n^2$ edges. 
Since each edge that was removed in this way was contained in at most $n^3$ copies of $C_5$, $\nu(C_5,G') \geq \nu(C_5,G) - \varepsilon_In^5$. Therefore, 
\begin{enumerate}
    \item $d(C_5,G') > OPT_k(C_5) - 5!\cdot \varepsilon_I - \varepsilon_F$,
    \item $d(C_5,G') > OPT_{k-1}(C_5) + 5! \cdot \varepsilon_I$.
\end{enumerate}
Using the previous two inequalities, along with the fact that $G'$ is $\overline{P_3}$-free, we will now show that $G'$ must be a complete $k$-partite graph. 

\begin{claim}\label{Gpartite}
$G'$ is a complete $k$-partite graph.
\end{claim}
\begin{proof}
Since $G'$ does not contain any induced copies of $\overline{P_3}$ as a subgraph, each pair of non-adjacent vertices must have an identical neighborhood. Therefore, we can partition $V(G')$ into independent sets $X_1,\ldots, X_{\ell}$ such that for all distinct $i,j \in [\ell]$, each vertex in $X_i$ is adjacent to each vertex in $X_j$. Hence, $G'$ is a complete $\ell$-partite graph. Since 
\[OPT_{k-1}(C_5)  = \lim\limits_{n \to \infty} d(C_5,T_{k-1}(n))\]
and $d(C_5,G') > OPT_{k-1}(C_5) + 5!\cdot\varepsilon_I$,  Proposition~\ref{setsize} implies that $G'$ must be $k$-partite if $n$ is sufficiently large.
\end{proof}

\medskip
At this point, we know that $G'$ only differs from $T_k(n)$ in the sizes of the partite sets $X_1,X_2,\ldots,X_k$. The next claim will show that we can impose that the partite sets in $G'$ must be reasonably close to being balanced. 
\begin{claim}\label{balanced partite sets}
Let $G'$ be a complete $k$-partite graph with partite sets $X_1,X_2,\dots,X_k$. Then for any $\varepsilon_T > 0$, there exists $\delta > 0$ such that if 
\[d(C_5,G') > OPT_k(C_5) - \delta,\]
then for each $i = 1,\dots,k$
\[
\frac{n(1 - \varepsilon_T)}{k} \leq |X_i| \leq \frac{n(1 + \varepsilon_T)}{k}.
\]
\end{claim}

\noindent 
\begin{proof}
For each $i = 1,\dots,k$ let $x_i = |X_i|$. 
Let $\varepsilon'(k-1) > \varepsilon_T$ and assume by symmetry that $x_1 = \frac{1+\varepsilon'(k-1)}{k}n$. If we picked $x_1 = \frac{1+\varepsilon'}{k}n$, we would get less pleasant expressions in what follows.
We want to calculate an upper bound on $d(C_5,G')$.
By Lemma~\ref{move a vertex}, $d(C_5,G')$ is maximized if all remaining parts are balanced. That is, $x_i = \frac{1-\varepsilon'}{k}n$ for $i = 2,\ldots,k$. 
With knowing the sizes of all $X_i$, the following is a straightforward calculation,
\begin{align*}
d(C_5,G') \leq &~
OPT_k(C_5) 
- 60\varepsilon'^2 \left(1-\frac6k + \frac{15}{k^2} - \frac{18}{k^3} + \frac{8}{k^4}\right) 
+ 60\varepsilon'^3 \left(1 - \frac{8}{k} + \frac{25}{k^2} - \frac{34}{k^3} + \frac{16}{k^4} \right) 
\\
&~
+ 180\varepsilon'^4 \left(\frac{1}{k} - \frac{5}{k^2} + \frac{8}{k^3} - \frac{4}{k^4}\right) 
- 12\varepsilon'^5 \left(1-\frac{15}{k^2} + \frac{30}{k^3} - \frac{16}{k^4} \right) + o(1), 
\end{align*}
see Appendix~\ref{appXX} for a code in SageMath. 

For all $k \geq 3$, the term $1-\frac6k + \frac{15}{k^2} - \frac{18}{k^3} + \frac{8}{k^4}$ is positive with minimum $\frac{8}{81}$ at $k=3$. 
For sufficiently small $\varepsilon'$ and large $n$, we get
\begin{align*}
d(C_5,G') \leq OPT_k(C_5) - 5 \varepsilon'^2.  
\end{align*}
This implies the statement of the claim.
\end{proof}

Returning to the proof of Lemma~\ref{friststab}, suppose that $\varepsilon > 0$, and let $\varepsilon_T = \varepsilon/2$. Next, choose an $\varepsilon_I \leq \varepsilon/2$ small enough so that $\varepsilon_F$ and $\delta_F$ are sufficiently small. In particular, we must select $\varepsilon_I$, $\varepsilon_F$, and $\delta_F$ so that any $k$-partite graph $G'$ satisfying $d(C_5,G') > OPT_k(C_5) - 5!\varepsilon_I - \varepsilon_F,$ must have partite sets $X_1,\dots,X_k$ that satisfy  
$$ \frac{n(1 - \varepsilon/2)}{k} \leq |X_i| \leq \frac{n(1 + \varepsilon/2)}{k}$$
for all $i = 1, \ldots, k$. Then by changing at most $(\varepsilon_I + \varepsilon_T)n^2$ pairs we can obtain $T_k(n)$ from the original graph $G$, which completes the proof of Lemma~\ref{friststab}.
\end{proof}

\section{Exact Result}
In this section we will prove Theorem~\ref{main}(ii). First we will give a brief outline. As we have shown, if $G$ is a $K_{k+1}$-free graph on $n$ vertices for large enough $n$ that contains close to the extremal number of five-cycles, then the edit distance between $G$ and $T_k(n)$ is very small. Given such a graph $G$, the process of deleting and adding the necessary edges to transform $G$ into the Tur\'{a}n graph actually increases the number of five-cycles. This will prove that $T_k(n)$ is the unique extremal graph for large enough $n$.

\begin{proof}[Proof of Theorem~\ref{main}(ii)]
Suppose that $k \geq 3$. By Lemma~\ref{friststab}, there exists an $\varepsilon > 0$ and an integer $n_0 = n(k,\varepsilon)$ so that for any $K_{k+1}$-free graph $G$ on $n \geq n_0$ vertices satisfying 
\[
d(C_5,G) > OPT_k(C_5) - \varepsilon,
\]
we have that $\Dist(G,T_k(n)) \leq \frac{1}{2k^{10}}n^2$. 

Let $G$ be a graph on $n$ vertices, where $n$ is sufficiently large. In particular, $n \geq \max\{n(k,\varepsilon),2k^5+1\}$ and $G$ satisfies
%
\[
d(C_5,G) > OPT_k(C_5) - \varepsilon',
\]
where $\varepsilon' \leq \min\{\varepsilon,\frac{1}{k^{10}}\}$.
Lemma~\ref{friststab} gives a partition of $V(G)$ into $k$ sets 
$X_1,X_2,\dots,X_k$, where $\lfloor \frac{n}{k} \rfloor \leq |X_i| \leq \lceil \frac{n}{k} \rceil$ for all $i = 1,\dots,k$, so that by changing at most $\frac{1}{2k^{10}} n^2$ pairs $uv$ for $u,v \in V(G)$, we can construct a new graph $G'$ from $G$ so that $G'$ is isomorphic to $T_k(n)$ and the partite sets of $G'$ are $X_1,X_2,\dots,X_k$. 

Call each edge that is removed in this process a \emph{surplus edge} and call each edge that is added in this process a \emph{missing edge}. For each vertex $v \in V(G)$, let $f_v$ denote the sum of the total number of surplus edges and missing edges incident to $v$. Define the set $X_0$ to contain each vertex $v$ with $f_v > \frac{1}{k^{6}} n$. We will refer to each vertex in $X_0$ as a \emph{bad vertex}.

\begin{claim}\label{xzerosize}
$|X_0| \leq \frac{1}{k^4} n$. 
\end{claim}

\begin{proof}[Proof]
Since $f_v > \frac{1}{k^6} n$ for each vertex $v \in X_0$ and the combined total of surplus edges and missing edges in $G$ is at most $\frac{1}{2k^{10}} n^2$, it follows that
$$\frac{1}{k^6}n|X_0| \leq \frac{1}{k^{10}}n^2,$$
which proves Claim~\ref{xzerosize}. 
\end{proof}

\medskip
For a graph $G$ and a vertex $v \in V(G)$ let $N_G(v)$ denote the neighborhood of $v$ in $G$. For all $v \in V(G)$, let $d_i(v)$ denote the size of the set $N_G(v) \cap (X_i \setminus X_0)$. Let 
$$d^*(v)  = \sum\limits_{i = 1}^k d_i(v).$$
By Claim~\ref{xzerosize}, 
$$\left\lfloor \frac{n}{k}\right\rfloor - \frac{1}{k^4} n \leq |X_i \setminus X_0| \leq \left\lceil \frac{n}{k} \right\rceil,$$
for all $i = 1,\ldots,k$.
Thus, for each vertex $v$ not contained in $X_0$,
$$d^*(v) \geq \left(\frac{k-1}{k} - \frac{1}{k^4} - \frac{1}{k^6}\right)n.$$
For two vertices $u$ and $v$ in a graph $G$, let $N_G(u,v)$ denote the \emph{common neighborhood} of $u$ and $v$, which is the set of all vertices in $G$ adjacent to both $u$ and $v$.
\begin{claim}\label{degup}
There are no surplus edges in $G - X_0$. 
\end{claim}

\begin{proof}
Assume by way of contradiction that $G - X_0$ contains a surplus edge $uv$.
Our goal is to show that it would be in $K_{k+1}$.
Since $uv$ is removed in the the process of transforming $G$ into the Tur\'{a}n graph, we may assume by symmetry that $u$ and $v$ are contained in the same set $X_1$. 
Since neither vertex is contained in $X_0$, 
\[
\min\{d_j(v),d_j(u)\} \geq \left(\frac{1}{k} - \frac{1}{k^4} -  \frac{1}{k^6}\right) n
\]
for each $j = 2,\dots,k$. Therefore,
$$|N_G(u,v) \cap (X_2 \setminus X_0)| \geq \left(\frac{1}{k} - \frac{1}{k^4} - \frac{2}{k^6}\right)n > 0.$$
Pick one vertex $w_2$ contained in $N_G(u,v) \cap (X_2 \setminus X_0)$. Since $w_2$ is not contained in $X_0$,
$$|N_G(w_2) \cap N_G(u,v) \cap (X_3 \setminus X_0)| \geq \left(\frac{1}{k} - \frac{1}{k^4} - \frac{3}{k^6}\right)n > 0.$$

This implies that we can find some common neighbor, say $w_3$, of $u,v,$ and $w_2$, where $w_3 \in X_3\setminus X_0$. 
We continue the process of a selecting a vertex $w_{j} \in X_j \setminus X_0$ in the common neighborhood of the set $\{u,v,w_1,\dots,w_{j-1}\}$ for all $j = 4,\dots,k$. This is possible because after selecting $w_{j-1}$, the common neighborhood of the set $\{u,v,w_2,\dots,w_{j-1}\}$ contains at least 
$$ \left(\frac{1}{k} - \frac{1}{k^4} - \frac{j}{k^6}\right)n > 0$$
vertices in $X_{j}\setminus X_0$ for all $j = 4,\dots,k$. 
This implies, however,  that the set $\{
u,v,w_2,\dots,w_{k}\}$ obtained by selecting a vertex in this way from each partite set $X_2,\dots,X_k$ induces a copy of $K_{k+1}$ in $G$, which is a contradiction. 
\end{proof}

\medskip
An immediate consequence of Claim~\ref{degup} is that every surplus edge in $G$ is incident to at least one vertex in $X_0$, implying that $G - X_0$ is a $k$-partite graph, albeit not necessarily complete $k$-partite. We will split the vertices of $X_0$ into two classes. For each vertex $v \in X_0$, one of the following holds. 
\begin{enumerate}
    \item There exists some index $i \in \{1,2,\dots,k\}$ such that $d_i(v) = 0$. In this case we will call $v$ a \emph{type 1} vertex, or
    \item $d_i(v) > 0$ for all $i = 1,\dots,n$. In this case we will call $v$ a \emph{type 2} vertex. 
\end{enumerate}

As we are trying to show that every extremal graph is a complete balanced $k$-partite graph, we will now prove that $G$ cannot contain any type $2$ vertices. First in Claim~\ref{upgooddeg}, we will prove that if $v$ is a type $2$ vertex, then $d^*(v)$ must be relatively small. In Claim~\ref{finishub}, we will prove a lower bound on the number of five-cycles containing a vertex $v$. Finally, in Claim~\ref{notype2}, we will show that a type 2 vertex cannot be contained in enough five-cycles to justify the claim that $G$ is an extremal graph. 
\begin{claim}\label{upgooddeg}
Let $v \in X_0$ be a type $2$ vertex. Then there exist distinct integers $i$ and $j$ where $1 \leq i,j \leq k$ such that
$$1 \leq d_{i}(v) \leq d_{j}(v) \leq \frac{1}{k^5}n.$$
\end{claim}

\begin{proof}
By symmetry, assume that $1 \leq d_1(v)$ and $d_1(v) \leq d_q(v) $ for all $q = 2,\dots,k$.
For contradiction, assume 
$d_q(v) > \frac{1}{k^5}n$ for all $q = 2,\dots,k$.
Let $w_1 \in X_1 \setminus X_0$ be adjacent to $v$ in $G$. Since $w_1 \notin X_0$,
\[
|N_G(v,w_1) \cap (X_2 \setminus X_0)| \geq \frac{1}{k^5}n - \frac{1}{k^6}n,
\]
implying that there exists a vertex $w_2 \in X_2\setminus X_0$ for which the set $\{v,w_1,w_2\}$ induces a triangle in $G$. 
If we continue selecting vertices in this way, then for all $q = 3,\dots,k$, there are at least 
$$\frac{1}{k^5}n - \frac{q-1}{k^6}n > 0$$
vertices in $X_q \setminus X_0$ that are adjacent to all of the previously selected vertices $v,w_1,\dots,w_{q-1}$. 
This implies that we can select $k$ vertices $w_1,\dots,w_k$ so that the set $\{v,w_1,\dots,w_k\}$ induces a copy of $K_{k+1}$ in $G$, which is a contradiction. Therefore, there exists an index $j \in \{2,\dots,k\}$ for which $1 \leq d_1(v) \leq d_j(v) \leq \frac{1}{k^5}n$, completing the proof of Claim~\ref{upgooddeg}.
\end{proof}

\begin{claim}\label{finishub} For all $k \geq 3$, and $v \in V(G)$, $\nu_G(v,C_5) \geq (OPT_k(C_5)-\frac{1}{k^{10}}) \binom{n}{4} - \frac{1}{k^5}n^4$.
\end{claim}

\begin{proof}
Suppose by way of contradiction that there exists some vertex $v$ for which 

$$\nu_G(v,C_5) < \left(OPT_k(C_5)-\frac{1}{k^{10}}\right) \binom{n}{4} - \frac{1}{k^5}n^4.$$

Since $d(C_5,G) > \OPT_k(C_5) - \frac{1}{k^{10}}$, it follows by averaging that there exists some vertex $u \in V(G)$ for which 
$$\nu_G(u,C_5) \geq \left(OPT_k(C_5)-\frac{1}{k^{10}}\right) \binom{n}{4}.$$ 
Let $\nu_G(\{u,v\},C_5)$ denote the number of five-cycles containing both $u$ and $v$. Then
$$\nu_G(\{u,v\},C_5) \leq 2n^3.$$
Let $G'$ be the graph obtained from $G$ by deleting $v$ and replacing it with a copy $u'$ of $u$. Since there is no edge between $u'$ and $u$, $G'$ is also $K_{k+1}$-free. As there were previously $\nu(\{u,v\},C_5)$ five-cycles containing $u$ and $v$,
\[
\nu(C_5,G') - \nu(C_5,G) \geq \nu_G(u,C_5) - \nu_G(v,C_5) - \nu_G(\{u,v\},C_5) \geq \frac{1}{k^5}n^4 - 2n^3 > 0
\]
since $n > 2k^5$. This, however, contradicts the assumption that $G$ is an extremal graph as $\nu(C_5,G') > \nu(C_5,G)$.
Therefore, if $n$ is sufficiently large it follows that for each $v \in V(G)$, 
\[\nu_G(v,C_5) \geq \left(OPT_k(C_5)-\frac{1}{k^{10}}\right) \binom{n}{4} - \frac{1}{k^5}n^4,
\]
which completes the proof of Claim~\ref{finishub}. 
\end{proof}

In Claims~\ref{notype2} and \ref{notype1} we count the number of 5-cycles containing a particular vertex $v \in X_0$. We use the following argument repeatedly.
We want to count the number of 5-cycles $vu_1u_2u_3u_4v$, where $v$ is in $X_0$, $u_1 \in X_i$, $u_4 \in X_j$ and $u_2,u_3 \in V(G) \setminus X_0$. 
 Assume we already picked $u_1$ and $u_4$ and want to count the number of choices for $u_2$ and $u_3$.
We distinguish two cases.
\begin{enumerate}
\item $i=j$: First $u_2$ can be in any of the remaining $k-1$ parts. Then $u_3$ has $k-2$ choices for a part to complete the $5$-cycle as it needs to avoid the parts containing $u_2$ and $u_4$ and these are distinct. After multiplying by $n^2/k^2$, the number of choices for $u_2$ and $u_3$ in each of the selected parts, we get
\[
\frac{(k-1)(k-2)}{k^2}n^2
\]
choices for $u_2$ and $u_3$ together.
\item $i\neq j$:
We further distinguish two cases. 
If $u_2 \notin X_j$, then there are $k-2$ parts which could contain $u_2$ and $k-2$ parts which could contain $u_3$. 
If $u_2 \in X_j$, then there are $k-1$ parts which could contain $u_3$.
After including the number of choices in each part, we get 
\[
\left( \frac{(k-2)^2}{k^2} + \frac{k-1}{k^2}\right)n^2
\]
choices  for $u_2$ and $u_3$ together.

\end{enumerate}

\begin{claim}\label{notype2}
$G$ does not contain any type 2 vertices. 
\end{claim}

\begin{proof}
Assume for contradiction that $v \in X_0$ is a type 2 vertex. Then by Claim~\ref{upgooddeg} there are two sets, say $X_1$ and $X_2$, such that 
$$1 \leq d_1(v) \leq d_2(v) \leq \frac{1}{k^5}n.$$
We will now provide an upper bound on the value of $\nu_G(v,C_5)$. We will count the maximum number of such five-cycles of the form $vu_1u_2u_3u_4v$ based on the locations of $u_1$ and $u_4$ as follows:
\begin{enumerate}
    \item If $u_1,u_4 \in X_1 \setminus X_0$ or $u_1,u_4 \in X_2 \setminus X_0$:
    \begin{equation}\label{type21}
    2\binom{\frac{n}{k^5}}{2}\frac{(k-1)(k-2)}{k^2}n^2.
    \end{equation}
    \item  $u_1 \in X_1 \setminus X_0$ and $u_4 \in X_2 \setminus X_0$:
    \begin{equation}\label{type22}
    \left(\frac{n}{k^5}\right)^2  \left( \frac{(k-2)^2}{k^2} + \frac{k-1}{k^2} \right)n^2.
    \end{equation}
    \item $u_1 \in (X_1 \setminus X_0) \cup  (X_2 \setminus X_0)$ and $u_4 \notin X_1 \cup X_2$:
    \begin{equation}\label{type23}
    \frac{2n}{k^5} \cdot \frac{n}{k} \left( \frac{(k-2)^2}{k^2} + \frac{k-1}{k^2} \right)n^2.
    \end{equation}
    \item $u_1,u_4 \notin X_1 \cup X_2$:
    \begin{equation}\label{type24}
    (k-2) \cdot \binom{\frac{n}{k}}{2} \cdot \frac{(k-1)(k-2)}{k^2}n^2 + \binom{k-2}{2} \cdot \frac{n^2}{k^2} \cdot \left( \frac{(k-2)^2}{k^2} + \frac{k-1}{k^2} \right)n^2.
    \end{equation}
\end{enumerate}
Finally, there are at most $\frac{2n^4}{k^4}$ five-cycles containing $v$ and at least one other vertex in $X_0$. Combining this, along with the upper bounds obtained in equations~\eqref{type21}--\eqref{type24},
\begin{equation*}
    \nu_G(v,C_5) \leq \frac{n^4}{24}\left(12 - \frac{84}{k} + \frac{228}{k^2} - \frac{300}{k^3} + \frac{216}{k^4} + \frac{48}{k^{6}} - \frac{144}{k^7} + \frac{144}{k^8} + \frac{48}{k^{10}} - \frac{144}{k^{11}} + \frac{120}{k^{12}}\right).
\end{equation*}
The SageMath code for verifying this fact can be found in Appendix~\ref{appendix code}. This implies that for large enough $n$, 
$$\left(OPT_k(C_5)-\frac{1}{k^{10}}\right)\binom{n}{4} - \nu_G(v,C_5) \geq \frac{1}{k^5}n^4.$$
Using SageMath, we verified that this was true for $3 \leq k \leq 1000$. After that, it is straightforward to check the coefficients in order to verify this fact. 
This contradicts Claim~\ref{finishub} since $G$ was assumed to be an extremal graph. Therefore, $G$ does not contain any type $2$ vertices. 
\end{proof}

\medskip
Since $G$ does not contain any type 2 vertices, we can place each vertex $v \in X_0$ into the set $X_i$ for which $d_i(v) = 0$. In order to show that $G$ is a complete $k$-partite graph, we must show that any pair of vertices $u$ and $v$ that were in $X_0$ and go to the same $X_i$  cannot be adjacent. The next claim will provide an upper bound on the ``good degree" of at least one of these adjacent vertices. 

\begin{claim}\label{type1deg}
Suppose that $u$ and $v$ are two adjacent type $1$ vertices such that $d_j(u) = d_j(v) = 0$ for some index $j \in \{1,\dots,k\}$. Then there exists some index $i \in \{1,\dots,k\}$ such that $i \neq j$ and
\[
d_i(u) \leq \frac{k^2 + 1}{2k^3}n \text{ or }  d_i(v) \leq \frac{k^2 + 1}{2k^3}n.
\]
\end{claim}
\begin{proof}
By symmetry we may assume that $j = 1$. Assume for contradiction that 
\[
| N_G(u,v) \cap (X_i\setminus X_0)| > \frac{1}{k^3}n
\]
for all $i = 2,\dots,k$. Using an identical argument to the one made in the proof of Claim~\ref{upgooddeg}, there exists a set $\{w_2,\dots,w_k\}$ such that $w_i \in (X_i \setminus X_0)$ and the set $\{u,v,w_2,\dots,w_k\}$ induces a $K_{k+1}$ in $G$, which is a contradiction. This implies that for at least one index $i$, 
\[
| N_G(u,v) \cap (X_i\setminus X_0)| \leq \frac{1}{k^3}n.
\]
Without loss of generality assume that $d_i(u) \leq d_i(v)$. Then
\[
d_i(u) \leq \frac{n}{2}\left(\frac{1}{k} - \frac{1}{k^3}\right) + \frac{n}{k^3}=\frac{k^2+1}{2k^3}n,
\]
which completes the proof of Claim~\ref{type1deg}. 
\end{proof}

\medskip
We will now show that the vertex $u$ of low degree described in the previous claim cannot be contained in enough five-cycles to justify the assumption that $G$ is an extremal graph. Unlike Claim~\ref{notype2}, we will only show that the two vertices $u$ and $v$ from Claim~\ref{type1deg} cannot be adjacent.

\begin{claim}\label{notype1}
Suppose that $u$ and $v$ are type $1$ vertices such that $d_j(u) = d_j(v) = 0$ for some $j = 1,\dots,k$. Then $u$ and $v$ are not adjacent. 
\end{claim}
\begin{proof}
By symmetry we may assume that $d_1(u) = d_1(v) = 0$. Assume for contradiction that $u$ and $v$ are adjacent. 
By symmetry and Claim~\ref{type1deg}, we may assume that $d_1(u) = 0$ and
\[
d_2(u) \leq \frac{k^2 + 1}{2k^3}n.
\]
In a similar manner as in Claim~\ref{notype2}, we will count the number of five-cycles of the form $uv_1v_2v_3v_4u$ incident to $u$ by considering the possibilities for the locations of $v_1$ and $v_4$ as follows:
\begin{enumerate}
    \item $v_1,v_4 \in X_2 \setminus X_0$:
    \begin{equation}\label{type11}
    \binom{\frac{k^2+1}{2k^3}n}{2}\frac{(k-1)(k-2)}{k^2}n^2.
    \end{equation}
    
    \item $v_1 \in X_2 \setminus X_0$ and $v_4 \notin X_2$:
    \begin{equation}\label{type12}
    \frac{k^2 + 1}{2k^3} \cdot \frac{k-2}{k}  \left( \frac{(k-2)^2}{k^2} + \frac{k-1}{k^2} \right)n^4.
    \end{equation}
    
    \item $v_1,v_4 \notin X_2$:
    \begin{equation}\label{type13}
    (k-2) \cdot \binom{\frac{n}{k}}{2} \cdot \frac{(k-1)(k-2)}{k^2}n^2 + \binom{k-2}{2} \cdot \frac{n^2}{k^2} \cdot \left( \frac{(k-2)^2}{k^2} + \frac{k-1}{k^2} \right)n^2.
    \end{equation}
\end{enumerate}
There are at most $\frac{2}{k^4}n^4$ five-cycles containing $u$ and at least one other vertex in $X_0$. 
Combining this along with equations~\eqref{type11}--\eqref{type13},
\[
\nu_G(u,C_5) \leq \left(12 - \frac{72}{k} + \frac{171}{k^2} - \frac{189}{k^3} + \frac{96}{k^4} + \frac{90}{k^5} + \frac{57}{k^6} - \frac{9}{k^7} + \frac{6}{k^8}\right)\frac{n^4}{24}.
\]
%
For $k > 1000$ it is clear that 
\[
12 - \frac{72}{k} + \frac{171}{k^2} - \frac{189}{k^3} + \frac{96}{k^4} + \frac{90}{k^5} + \frac{57}{k^6} - \frac{9}{k^7} + \frac{6}{k^8} \leq \OPT_k(C_5) - \frac{1}{k^{10}}.
\]
The SageMath code for verifying that this is also true for $3 \leq k \leq 1000$ found in Appendix~\ref{appendix code type 1}.
Given this fact, it is straightforward to verify that
$$\left(\OPT_k(C_5) - \frac{1}{k^{10}}\right)\binom{n}{4} - \nu_G(u,C_5) \geq \frac{1}{k^5}n^4$$
for large enough $n$. This, however contradicts Claim~\ref{finishub}, which implies that $u$ and $v$ are not adjacent. 
\end{proof}

\medskip
Claim~\ref{notype1} implies that if $u$ and $v$ are type $1$ vertices for which $d_i(v) = d_i(u) = 0$, then $u$ and $v$ cannot be adjacent. This means that we can place each type $1$ vertex $v$ into the set $X_i$ for which $d_i = 0$. Since $G$ does not contain any type $2$ vertices, this implies that $G$ is a $k$-partite graph. Since $G$ maximizes the number of (possibly non-induced) $C_5$ subgraphs, it follows that $G$ must be a compelte $k$-partite graph. Finally, Proposition~\ref{setsize} implies that $G$ is isomorphic to $T_k(n)$, implying that for large enough $n$, the Tur\'{an} graph $T_k(n)$ is the unique extremal graph maximizing the number of $C_5$ subgraphs. \end{proof}

\section{Conclusion}
In \cite{palmer}, Palmer and Gerbner showed that if $H$ is a graph and $F$ is a graph with chromatic number $k+1$, then 
$$\ext(n,H,F) \leq \ext(n,H,K_{k+1}) + o(n^{|H|}).$$
Since the Tur\'{a}n graph $T_k(n)$ does not contain any $(k+1)$-chromatic graph as a subgraph, this immediately implies that for any $(k+1)$-chromatic graph $F$,
$$\lim\limits_{n \to \infty}d(C_5,F) = \frac{1}{k^4}(12k^4 - 60k^3 + 120k^2 - 120k + 48),$$
which closely resembles the Erd\H{o}s-Stone-Simonovits theorem. 

Let $G$ be a graph with chromatic number $k$. Then for any $r \geq k$, the Tur\'{a}n graph $T_r(n)$ contains $G$ as a subgraph. When trying to maximize the copies of $G$ among $K_{r+1}$-free graphs, evidence seems to suggest that $T_r(n)$ is extremal, as we have shown to be the case with five-cycles. While a complete $r$-partite graph seems to frequently be the best option, it is not always optimal to balance the partite sets.  
Let $S_t$ be a star with $t$ leaves, also known as $K_{1,t}$.
Is it easy to see $\ext(n,S_4,K_3)$ is achieved by an unbalanced bipartite graph.
For more detailed treatment of stars, see  Cutler, Nir, and Radcliffe~\cite{stars}.
%
%
%
It seems very likely that while the Tur\'{a}n graph is not always extremal, that some complete $r$-partite graph will be best possible. 

\begin{conj}
Let $G$ be a graph and let $k > \chi(G)$ be an integer. Then for all $r \geq k$, $\ext(n,G,K_r)$ is realized by a complete $(r-1)$-partite graph. 
\end{conj}

While an unbalanced $r$-partite graph might best possible in some cases, we believe that for large enough $r \geq \chi(G)$, the value of $\ext(n,G,K_{r+1})$ is realized by the Tur\'{a}n graph. As $r$ increases, any $G$-subgraph in $T_r(n)$ can be taken from an increasing number of partite sets. Thus, as $r$ grows larger, the effect of $G$ being unbalanced becomes minimized. The following conjecture also appears in~\cite{gerbner2020exact}.

\begin{conj}
Let $G$ be a graph and let $r > |V(G)|$ be an integer. Then $\ext(n,G,K_r)$ is realized by the Tur\'{a}n graph $T_{r-1}(n)$.  
\end{conj}

\section*{Acknowledgement}
We would like to thank the anonymous referees for their comments that greatly improved the presentation of the paper.

\bibliographystyle{abbrv}
\bibliography{C5Krfree}
\newpage
\section{Appendix}
SageMath code on next pages can be also obtained at \url{https://arxiv.org/abs/2007.03064}.

\begin{table}[h]
\begin{center}
\medskip
$F_0 = \Fa$ \hspace{2 em} $F_1 = \Fb$ \hspace{2 em} $F_2 = \Fc$ \hspace{2 em} $F_3 = \Fd$ 

\medskip
$F_4 = \Fe$ \hspace{2 em} $F_5 = \Ff$ \hspace{2 em} $F_6 = \Fg$ \hspace{2 em} $F_7 = \Fh$ 

\medskip
$F_8 = \Fi$ \hspace{2 em} $F_9 = \Fj$ \hspace{2 em} $F_{10} = \Fk$ \hspace{2 em} $F_{11} = \Fl$ 

\medskip
$F_{12} = \Fm$ \hspace{2 em} $F_{13} = \Fn$ \hspace{2 em} $F_{14} = \Fo$ \hspace{2 em} $F_{15} = \Fp$

\medskip
$F_{16} = \Fq$ \hspace{2 em} $F_{17} = \Fr$ \hspace{2 em} $F_{18} = \Fs$ \hspace{2 em} $F_{19} = \Ft$

\medskip
$F_{20} = \Fu$ \hspace{2 em} $F_{21} = \Fv$ \hspace{2 em} $F_{22} = \Fw$ \hspace{2 em} $F_{23} = \Fx$

\medskip
$F_{24} = \Fy$ \hspace{2 em} $F_{25} = \Fz$ \hspace{2 em} $F_{26} = \Faa$ \hspace{2 em} $F_{27} = \Fbb$

\medskip
$F_{28} = \Fcc$ \hspace{2 em} $F_{29} = \Fdd$ \hspace{2 em} $F_{30} = \Fee$ \hspace{2 em} $F_{31} = \Fff$

\medskip
$F_{32} = \Fgg$ \hspace{2 em} $F_{33} = \Fhh$
\end{center}
\caption{Graphs on 5 vertices up to isomorphism.}
\label{tab:fivevertexgrahs}
\end{table}

\begin{table}[h]

$H_{4,1}^{\sigma_1} = \fourone$ \hspace{2 cm} $H_{4,2}^{\sigma_1} = \fourtwo$

$H_{4,3}^{\sigma_2} = \fourthree$ \hspace{2 cm} $H_{4,4}^{\sigma_2} = \fourfour$ \hspace{2 cm} $H_{4,5}^{\sigma_2} = \fourfive$

$H_{4,6}^{\sigma_3} = \foursix$ \hspace{2 cm}
$H_{4,7}^{\sigma_3} = \fourseven$ \hspace{2 cm}
$H_{4,8}^{\sigma_3} = \foureight$

\caption{Labeled graphs on four vertices.}\label{tab:fourvertex}
\end{table}

\newpage
\subsection{Proof of Claim~\ref{claim that C1 is max}}\label{proof that C1 is max}

\begin{lstlisting}[language=Python]
# SageMath code for Claim 2.1

var('k')

# Vector containing coefficients for each graph in F_5
cFi = [0]*34

# Z(K)
Zk = [(k*k*k*k-10*k*k*k+35*k*k-50*k+24)/k^4]*34
Zk[33] = (-10*k*k*k+35*k*k-50*k+24)/k^4

# P_i(K) first init with 0
P1k = [0]*34
P2k = [0]*34
P3k = [0]*34
P4k = [0]*34
P5k = [0]*34

# P_1(K)
P1k[0] = 10*k^2 - 20*k + 10
P1k[1] = k^2 -2*k +1
P1k[3] =  -k + 1
P1k[4] = -4*k + 4
P1k[18] = 1
P1k[19] = 1

# P_2(K)
P2k[18] = 3*k^2 - 12*k + 12
P2k[29] = k^2 - 6*k + 8
P2k[31] = -4*k+10
P2k[32] = 3

# P_3(K)
P3k[4] = 6*k^2 - 24*k + 24
P3k[11] = k^2 - 4*k + 4
P3k[17] = -k+2
P3k[19] = -6*k+12
P3k[31] = 2
P3k[32] = 3

# P_4(K)
P4k[19] =  6
P4k[28] = -1
P4k[30] =  2
P4k[31] = -4 

# P_5(K)
P5k[19] =  6*k^2-36*k+54
P5k[30] =  2*k^2-20*k+42
P5k[31] =  4*k^2 - 24*k + 36
P5k[32] =  -24*k + 84
P5k[33] =  120

## Scaling functions z(k) and p_i(k)
zk = 6*(5*k^3 - 20*k^2 + 30*k - 16)/(5*k^3 - 35*k^2 + 75*k - 48)
p1k = 3*(k^5 - 8*k^4 + 22*k^3 - 24*k^2 + 8*k)/ \
      (5*k^7 - 35*k^6 + 75*k^5 - 48*k^4)
p2k = (10*k^5 - 60*k^4 + 109*k^3 - 76*k^2 + 18*k)/ \
      (5*k^7 - 35*k^6 + 75*k^5 - 48*k^4)
p3k = (5*k^5 - 28*k^4 + 45*k^3 - 28*k^2 + 6*k)/ \
      (5*k^7 - 35*k^6 + 75*k^5 - 48*k^4)
p4k = (1/4)*(5*k^7 - 30*k^6 + 53*k^5 - 52*k^4 + 94*k^3 - 96*k^2 + 24*k)/ \
      (5*k^7 - 35*k^6 + 75*k^5 - 48*k^4)
p5k = (1/4)*(15*k^5 - 60*k^4 + 78*k^3 - 40*k^2 + 8*k)/ \
      (5*k^7 - 35*k^6 + 75*k^5 - 48*k^4)

# Number of C_5s  in each of the 34 graphs. 
five_cycles = [0]*26+[1,1,1,2,2,4,6,12]

def test_positive_for_small_k(x):
    for r in [4..1000]:
        if x.substitute(k=r) <= 0:
            print ("ERROR:",x,"is not positive for k=",r)
            return

# Denominator for all cFi in the result
den = 5*k^7 - 35*k^6 + 75*k^5 - 48*k^4

# Test if denominator is positive for small k
test_positive_for_small_k(den)

for i in [0..33]:
	cFi[i] = expand(factor( (zk*(Zk[i]) + p1k*P1k[i] + p2k*P2k[i] +\
      p3k*P3k[i] + p4k*P4k[i] + p5k*P5k[i] + five_cycles[i] )*den ))

# This is the polynomial for each of the tight graphs (when multiplied
# by den)
# It is the coefficient at c_{F_0} = cFi[0]
opt_pol = cFi[0]
#opt_pol = 60*k^7 - 720*k^6 + 3600*k^5 - 9876*k^4 + 16320*k^3 \
#          - 16440*k^2 + 9360*k - 2304

# Printing of all coefficients c_{F_i}
print("Printing the resulting coefficients")
for i in [0..33]:
	print ("cF{}= {} / {}".format(i,cFi[i],den))

# Test that  cFi[0] == optimum value
optimum = -60/k + 120/k^2 - 120/k^3 + 48/k^4 + 12

########### Start of Claim 2.1
# Calculating the differences C_1 - C_i for i=2,..,10.
# We do it by checking c_{F_0} - c_{F_i} for all i
#
cF0_cFi=[0]*34
for i in [0..33]:
	cF0_cFi[i] = expand(cFi[i]- opt_pol)

# Showing that c_F0 is largest by displaying the difference.
# The leading coefficient at k^7 is negative
# and tests for small values of k by evaluation

def test_not_positive_for_small_k(x):
    for r in [4..1000]:
        if x.substitute(k=r) > 0:
            print ("ERROR:",x,"is positive for k=",r)
            return
print()
print("Showing that c_F0 is largest")
print(" - see difference is 0 or leading coefficient negative")
for i in [0..33]:
    print ("numerator(c_F0 - c_F{})=".format(i),cF0_cFi[i])
    test_not_positive_for_small_k(cF0_cFi[i])


print("all done")
\end{lstlisting}

\newpage

\begin{lstlisting}[language=Python]

# SageMath code for Claim 2.1 when k=3

# This code calculates the coefficients cFi
#for the case of k = 3 in Claim 2.1.
# In \mathcal{F}_5: F20, F24, F30, F32, and F33 all contain a K4. 
# Here we have removed those graphs and re-indexed the remaining graphs.

cFi = [0]*29

#Counting the number of C5s in each graph.

constants = [0]*29
constants[24] = 1
constants[25] = 1
constants[26] = 1
constants[27] = 2
constants[28] = 4


#  Pi(K) first init with 0

P1k = [0]*34
P2k = [0]*34
P3k = [0]*34
P4k = [0]*34
P6k = [0]*34

P1k[0] = 40
P1k[1] = 4
P1k[3] =  -2
P1k[4] = -8
P1k[18] = 1
P1k[19] = 1

P2k[18] = 3
P2k[27] = -1
P2k[28] = -2

P3k[4] = 6
P3k[11] = 1
P3k[17] = -1
P3k[19] = -6
P3k[28] = 2

P4k[19] = 6
P4k[26] = -1
P4k[28] = -4 

P6k[11] = 1
P6k[23] = 2
P6k[22] = -1
P6k[27] = -2
P6k[17] = 1
P6k[19] = 6
P6k[28] = -4

#The scaling coefficients for k = 3

p1k = 1/27
p2k = 13/27
p3k = 8/27
p4k = 2/9
p6k = 17/54

# This calculates cFi
for i in [0..28]:
	cFi[i] = p1k*P1k[i] + p2k*P2k[i] + p3k*P3k[i] 
	\+ p4k*P4k[i] + p6k*P6k[i] + constants[i]

# This prints the value of cFi along with the (possibly) re-indexed graph.
for i in [0..28]:
	print ('coefficient of F',i,'=',cFi[i])

\end{lstlisting}

\newpage
\subsection{SageMath code for Claim~\ref{balanced partite sets}}\label{appXX}

\begin{lstlisting}[language=Python]
# SageMath code Claim 3.10

var('k,e')

# this is the size of the sets. 
# We start by using epsilon*(k-1) for easier counting. 

x = (1 + e*(k-1))/k
y = (1 - e)/k

# These count the number of five cycles. 
# The first is the one we use.
# The second is a sanity check. 
def fivecyclecount(x,y):

    # here we count by picking one vertex in x, 
    # then counting the number of possible five cycles.
    
    neighbors_in_same_sets = \
    x*(k-1)*(y^2/2)*(y^2*(k-2)*(k-3) + x*(k-2)*y*2)
    
    neighbors_in_diff_sets = \
    x*( (y*(k-1))*((k-2)*y) )/2*( (k-3)*(k-3)*y^2 
    + (k-2)*y^2 + x*(k-3)*y)
    
    # This is counting the number of five-cycles not in X_1.
    # Note that it is equal to the sanity check but with k-1.
    
    nobadset_twosame = (y^3*(k-1)*(k-2)/2)*(y^2*(k-2)*(k-3))
    nobadset_nosame =  \
    (y^3*(k-1)*(k-2)*(k-3)/2)*( (k-3)*(k-3)*y^2 
    + (k-2)*y^2) 
    
    return 120*(neighbors_in_diff_sets + \
    neighbors_in_same_sets) + \
    24*(nobadset_twosame + nobadset_nosame)

def sanity(y):
    nobadset_twosame = (y^3*k*(k-1)/2)*(y^2*(k-1)*(k-2))
    nobadset_nosame = \
    (y^3*k*(k-1)*(k-2)/2)*( (k-2)*(k-2)*y^2 + (k-1)*y^2)
    return 24*(nobadset_twosame + nobadset_nosame)
    
f = fivecyclecount(x,y)

# to check our count is correct, 
# notice that the non-epsilon terms equal OPT. 
view(f.collect(e))
view(expand(sanity(1/k)))
\end{lstlisting}

\newpage
\subsection{SageMath code for Claim 4.5}\label{appendix code}

\begin{lstlisting}[language=Python]
# Sage code for claim 4.5 - showing there are no type 2 vertices

var('k')

#This function counts the number of five cycles using equations (15) - (18) 

def fivecyclecount(k):
	onebadset_twobadvertices = 2*( 1/(2*k^10) )*( (k-1)*(k-2)/k^2 ) 
	twobadset_twobadvertices = \ 
	(1/k^10 + 2/k^6)*( (k-2)^2/k^2 + (k-1)/k^2 )
	nobadset_twosame = ( (k-2)/(2*k^2) )*( (k-1)*(k-2)/k^2 )
	nobadset_nosame = \ 
	( (k-2)*(k-3)/(2*k^2) )*( (k-2)^2/k^2 + (k-1)/k^2 )
	with_X0 = 2/k^4
	return 24*(onebadset_twobadvertices \ 
	+ twobadset_twobadvertices + nobadset_twosame \ 
	+ nobadset_nosame + with_X0)


# This gives the sum of equations (15) - (18) factored in a nice way. 

expanded_first_check = expand(fivecyclecount(k))
print('five cycles containing a type 2 vertex:\n',expanded_first_check)


# actual upper bound once we account for the vertices in X_0
def bad_ub(k):
	return  fivecyclecount(k)

# The average "density" of five cycles containing a particular vertex
def good_ub(k):
	return -60/k + 120/k^2 - 120/k^3 + 48/k^4 + 12 - 1/k^10

# The difference between the optimal and the count for type 2.
def epsilon(k):
	return good_ub(k) - bad_ub(k)

print('difference:',factor(epsilon(k)))

# This verifies that for small values of k, count(r) greater than 1/k^5 
def count_check(a,b):
	for i in [a..b]:
		if epsilon(i) < 1/(i^5):
			return "the difference is less than 1/k^5 for k =",i
	return "'difference' is greater than 1/k^5 for all values"+\
        " of "+str(a)+" <= k up to "+str(b)

print(count_check(3,1000))
\end{lstlisting}
\newpage
\subsection{SageMath code for Claim 4.7}\label{appendix code type 1}

\begin{lstlisting}[language=Python]
# SageMath code for Claim 4.7

var('k')

# This function counts the number of five-cycles in equations (19) - (21)
def fivecyclecount(k):
    onebadset_twobadvertices = \
    ( (1/2)*( (k^2+1)/(2*k^3) )^2 )*( (k-1)*(k-2)/k^2 )
    onebadset_onebadvertices = \
    ( (k^2+1)/(2*k^3) )*( (k-2)^3/k^3
    + (k-1)*(k-2)/k^3 )
    nobadset_twosame = \
    ( (k-2)/(2*k^2) )*( (k-1)*(k-2)/k^2 )
    nobadset_nosame = \
    ( (k-2)*(k-3)/(2*k^2) )*( (k-2)^2/k^2 + (k-1)/k^2 )
    with_X0 =  2/k^(4)
    return 24*(onebadset_twobadvertices \
    + onebadset_onebadvertices + nobadset_twosame \
    + nobadset_nosame+with_X0)


# This gives the sum of equations (19) - (21) factored in a nice way

expanded_first_check = expand(fivecyclecount(k))
print('five_cycles_containing_a_type1vertex:\n',expanded_first_check)

# The upper bound on five cycles containing a suboptimal type 1 vertex

def bad_ub(k):
    return  fivecyclecount(k)

def good_lb(k):
    return -60/k + 120/k^2 - 120/k^3 + 48/k^4 + 12 - 1/k^10

#difference from optimal value
def epsilon(k):
    return good_lb(k) - bad_ub(k)

print('difference:',factor(epsilon(k)))


# This verifies that for small values of k, epsilon(k) is greater than 1/k^5
def count_check(a,b):
    for i in [a..b]:
        if epsilon(i) < 1/(i^(5)):
            return "the difference is less than 1/k^5 for k =",i
    return "'difference' is greater than 1/k^5 for all values"+\
        " of "+str(a)+" <= k up to "+str(b)

print(count_check(3,1000))
\end{lstlisting}

\end{document}